\documentstyle{article}

\newread\epsffilein    % file to \read
\newif\ifepsffileok    % continue looking for the bounding box?
\newif\ifepsfbbfound   % success?
\newif\ifepsfverbose   % report what you're making?
\newif\ifepsfdraft     % use draft mode?
\newdimen\epsfxsize    % horizontal size after scaling
\newdimen\epsfysize    % vertical size after scaling
\newdimen\epsftsize    % horizontal size before scaling
\newdimen\epsfrsize    % vertical size before scaling
\newdimen\epsftmp      % register for arithmetic manipulation
\newdimen\pspoints     % conversion factor
\def\epsfbox#1{\global\def\epsfllx{72}\global\def\epsflly{72}%
    \global\def\epsfurx{540}\global\def\epsfury{720}%
    \def\lbracket{[}\def\testit{#1}\ifx\testit\lbracket
    \let\next=\epsfgetlitbb\else\let\next=\epsfnormal\fi\next{#1}}%
\def\epsfgetlitbb#1#2 #3 #4 #5]#6{\epsfgrab #2 #3 #4 #5 .\\%
    \epsfsetgraph{#6}}%
\def\epsfnormal#1{\epsfgetbb{#1}\epsfsetgraph{#1}}%
\def\epsfgetbb#1{%
\openin\epsffilein=#1
\ifeof\epsffilein\errmessage{I couldn't open #1, will ignore it}\else
    {\epsffileoktrue \chardef\other=12
     \def\do##1{\catcode`##1=\other}\dospecials \catcode`\ =10
     \loop
        \read\epsffilein to \epsffileline
        \ifeof\epsffilein\epsffileokfalse\else
           \expandafter\epsfaux\epsffileline:. \\%
        \fi
    \ifepsffileok\repeat
    \ifepsfbbfound\else
     \ifepsfverbose\message{No bounding box comment in #1; using defaults}\fi\fi
    }\closein\epsffilein\fi}%
\def\epsfclipoff{\def\epsfclipstring{\ifepsfdraft\space clip\fi}}%
\epsfclipoff
\def\epsfsetgraph#1{%
    \epsfrsize=\epsfury\pspoints
    \advance\epsfrsize by-\epsflly\pspoints
    \epsftsize=\epsfurx\pspoints
    \advance\epsftsize by-\epsfllx\pspoints
    \epsfxsize\epsfsize\epsftsize\epsfrsize
    \ifnum\epsfxsize=0 \ifnum\epsfysize=0
       \epsfxsize=\epsftsize \epsfysize=\epsfrsize
       \epsfrsize=0pt
      \else\epsftmp=\epsftsize \divide\epsftmp\epsfrsize
        \epsfxsize=\epsfysize \multiply\epsfxsize\epsftmp
        \multiply\epsftmp\epsfrsize \advance\epsftsize-\epsftmp
        \epsftmp=\epsfysize
        \loop \advance\epsftsize\epsftsize \divide\epsftmp 2
        \ifnum\epsftmp>0
           \ifnum\epsftsize<\epsfrsize\else
              \advance\epsftsize-\epsfrsize \advance\epsfxsize\epsftmp \fi
        \repeat
        \epsfrsize=0pt
      \fi
    \else \ifnum\epsfysize=0
      \epsftmp=\epsfrsize \divide\epsftmp\epsftsize
      \epsfysize=\epsfxsize \multiply\epsfysize\epsftmp
      \multiply\epsftmp\epsftsize \advance\epsfrsize-\epsftmp
      \epsftmp=\epsfxsize
      \loop \advance\epsfrsize\epsfrsize \divide\epsftmp 2
      \ifnum\epsftmp>0
         \ifnum\epsfrsize<\epsftsize\else
            \advance\epsfrsize-\epsftsize \advance\epsfysize\epsftmp \fi
      \repeat
      \epsfrsize=0pt
     \else
      \epsfrsize=\epsfysize
     \fi
    \fi
    \ifepsfverbose\message{#1: width=\the\epsfxsize, height=\the\epsfysize}\fi
    \epsftmp=10\epsfxsize \divide\epsftmp\pspoints
    \vbox to\epsfysize{\vfil\hbox to\epsfxsize{%
       \ifnum\epsfrsize=0\relax
         \includegraphics{\ifepsfdraft}%
       \else
         \epsfrsize=10\epsfysize \divide\epsfrsize\pspoints
         \includegraphics{\ifepsfdraft}%
       \fi
       \hfil}}%
\global\epsfxsize=0pt\global\epsfysize=0pt}%
{\catcode`\%=12 \global\let\epsfpercent=%\global\def\epsfbblit{%BoundingBox}}%
\long\def\epsfaux#1#2:#3\\{\ifx#1\epsfpercent
    \def\testit{#2}\ifx\testit\epsfbblit
       \epsfgrab #3 . . . \\%
       \epsffileokfalse
       \global\epsfbbfoundtrue
    \fi\else\ifx#1\par\else\epsffileokfalse\fi\fi}%
\def\epsfempty{}%
\def\epsfgrab #1 #2 #3 #4 #5\\{%
\global\def\epsfllx{#1}\ifx\epsfllx\epsfempty
       \epsfgrab #2 #3 #4 #5 .\\\else
    \global\def\epsflly{#2}%
    \global\def\epsfurx{#3}\global\def\epsfury{#4}\fi}%
\def\epsfsize#1#2{\epsfxsize}
\newtheorem{sect}{}[section]

\newcommand{\Hom}{\mbox{\rm Hom}}
\newcommand{\lk}{\mbox{$\ell k$}}
\begin{document}

\title{Quandle Cohomology and State-sum Invariants
of Knotted Curves and Surfaces \\
{\large Dedicated to Professor Kunio Murasugi for his $70$th birthday} }

\author{
J. Scott Carter
\\University of South Alabama \\
Mobile, AL 36688 \\ carter@mathstat.usouthal.edu
\and
Daniel Jelsovsky
\\ University of South Florida
\\ Tampa, FL 33620  \\ jelsovsk@math.usf.edu
\and
Seiichi Kamada
\\  Osaka City University \\
Osaka 558-8585, JAPAN\\
kamada@sci.osaka-cu.ac.jp
\and
Laurel Langford \\
University of Wisconsin at River Falls \\
River Falls, WI 54022 \\
laurel.langford@uwrf.edu
\and
Masahico Saito
\\ University of South Florida
\\ Tampa, FL 33620  \\ saito@math.usf.edu
}

\maketitle

\begin{abstract}
The 2-twist spun trefoil is an example
of a sphere that is knotted in 4-dimensional space.
A proof is given in this paper that this sphere is
distinct from the same sphere with
its orientation reversed.
Our proof is based on
a state-sum invariant for
knotted surfaces
developed
via a cohomology theory of racks and quandles
  (also known as distributive groupoids).

A quandle is a set with a binary operation ---
the axioms of which model the Reidemeister moves
in
classical knot
theory.
Colorings of diagrams of knotted curves and surfaces
by quandle elements, together with
cocycles of quandles, are used
to define state-sum invariants
for knotted circles in $3$-space and knotted surfaces in $4$-space.

Cohomology  groups of various quandles are
computed herein and applied to the study of
the state-sum invariants.
Non-triviality of the invariants are proved for variety of knots and links,
and conversely, knot invariants are used to prove
non-triviality of cohomology for
a variety of quandles.

\vspace{1in}



\end{abstract}

\newpage

\section{Introduction}

A quandle is
a  set  with  a self-distributive ($(a*b)*c = (a*c)*(b*c)$)
binary operation
the axioms of which are
partially
motivated by  classical knot theory.
We derive a cohomology theory for quandles
diagrammatically from Reidemeister moves for classical knots and
knotted
surfaces. Our definition
of quandle (co)homology is a modification of rack
(co)homology defined in  \cite{FRS1} and \cite{FRS2}.
Quandle cocycles are used to define state-sum invariants
for knots and links in dimension
  $3$ and for knotted surfaces in dimension $4$.
As the main application of the invariant,
we show that the invariant  detects non-invertible knotted surfaces.

The invariants defined are demonstrated to be non-trivial
on a variety of examples.
  In many cases, the
invariant is related to
  linking numbers
(Sections~\ref{danssect} and \ref{sss2}). In the case of a
3-component surface link,
there is a notion of $3$-fold linking, defined combinatorially,
that can be used to compute the invariant over trivial quandles
(Section~\ref{sss2}).
In the classical case of knotted curves,
it is shown that
the trefoil ($3_1$ in the tables)
and figure 8 knot ($4_1$ in the tables) have non-trivial
  (mod $2$)-cocycle invariants
over a $4$-element quandle associated to
the rotations
of a tetrahedron.
Conversely, knots are used to prove algebraic results
--- non-triviality of cohomology groups for a
variety of quandles.
As a main topological application, the 2-twist spun trefoil is shown
(Section~\ref{themain})
to be non-invertible, {\it i.e.},
  distinct from itself with
the reversed  orientation,  by evaluating the state-sum
invariant with a $3$-cocycle
  over the three element
dihedral quandle (defined below).

In  \cite{FRS1} and \cite{FRS2},
the general framework for defining invariants from racks and
quandles and their homology and cohomology
is outlined.
The present paper defines knot invariants by means of
a state-sum,
using quandle cocycles.
This cocycle invariant also can be seen as
an analogue of
  the Dijkgraaf-Witten
invariants for $3$-manifolds \cite{DW} in that colorings and
cocycles are used to define state-sum invariants.
Another
analogue
of the Dijkgraaf-Witten invariants
was applied
to triangulated $4$-manifolds
in \cite{CKS2}.
The non-invertibility for certain classical knots had long been presumed
since the  1920's
but proved
first by Trotter in the 1960's and later using hyperbolic
structures
(see \cite{Hart,Kawa,Kawa:book}).
Fox \cite{FoxTrip} presented a  non-invertible knotted sphere
using Alexander modules.
Alexander modules, however,
fail to  detect non-invertibility of
the $2$-twist spun
trefoil.
In this paper we show its non-invertibility
using the cocycle state-sum invariants.
In particular, the cocycle invariants are the first state-sum
invariants in dimension $4$ that  carry information
  not contained in the
Alexander modules.
It was pointed out
to us
by D. Ruberman that Farber-Levine pairings
\cite{Farber1975,Farber1977,Levine1977} 
and
  Casson-Gordon invariants
detect non-invertibility of some twist-spun knots \cite{Hill,Ruber}.
A. Kawauchi pointed out that our invariant
detects the non-invertibility of the twist-spun trefoil even after
adding trivial $1$-handles, thus increasing the genus of the surface.
Therefore, our invariant implies
the new topological results, that these higher genus surfaces
are non-invertible.   He also informed that
Farber-Levine parings are generalized for
knotted surfaces of higher genus \cite{Kawauchi1990a} (cf. \cite{Sekine1989a})
by use of his duality \cite{Kawauchi1986a},
and that the (Farber-Levine-Kawauchi)
parings  \cite{Kawauchi1990a} also detect non-invertibility of higher genus 
surfaces.
Thus relations between these invariants and the state-sum invariants deserve
investigation.

Our inspiration for the definition of these invariants is found
  in Neuchl's paper
\cite{Neu} where related cocycles are used to
give examples of  representations of a Hopf category
in a braided monoidal 2-category
using quantum groups of finite groups.
Our definition was derived from
an attempt to construct a $2$-functor from
the braided $2$-category of knotted surfaces
  as summarized in \cite{BL1}
and presented in detail in \cite{BL},
to another $2$-category constructed from quandles.

\begin{sect}{\bf Organization.\ }
{\rm
Section~\ref{Rack} contains the basic definitions
of racks and quandles.
Rack cohomology and quandle cohomology
are defined in Section~\ref{coho}.
Section~\ref{classical}  defines invariants of classical
knots and links via assigning 2-cocycles to crossings.
Section~\ref{surface}
contains the analogous definition  for knotted surfaces.
Section~\ref{quandleco}
presents  calculations of cohomology groups for some exemplary  quandles.
Section~\ref{reltogroup} relates the quandle
2-cocycles to group 2-cocycles when the
quandle is a group with conjugation as the operation.
Section~\ref{danssect} contains
computations in the case of classical knots and links.
Section~\ref{sss2} defines a notion of linking for knotted surfaces.
This linking is
used to exemplify non-triviality of
the state-sum invariant in the case of surfaces in 4-space.
Section~\ref{surfbraid} develops techniques for computation for surface braids.
In Section~\ref{themain} these techniques are applied to the 2-twist spun
  trefoil and its orientation reversed image to demonstrate that these
  knotted surfaces are distinct.
}\end{sect}

\begin{sect}{\bf Acknowledgements.\/} {\rm
We are grateful for a grant for visitors
  from Alabama EPSCoR's Mathematical Infrastructure Committee which
brought Masahico Saito and Laurel Langford to Mobile for discussions.
We have had valuable conversations with
J. Baez, J. Birman, R. Fenn, L. Kauffman, C. Rourke, B. Sanderson,
D. Ruberman, and D. Silver.
Jos\'{e} Barrionuevo,
Edwin Clark, and Cornelius Pillen
had helpful programming hints
  for the
computation of quandle cocycles.
Seiichi Kamada is being supported by a Fellowship from the
Japan Society for the Promotion
of Science.
}\end{sect}

\section{Racks, Quandles, and Knots}
\label{Rack}
A {\it quandle}, $X$, is a set with a binary operation $(a, b) \mapsto a * b$
such that

(I) For any $a \in X$,
$a* a =a$.

(II) For any $a,b \in X$, there is a unique $c \in X$ such that
$a= c*b$.

(III)
For any $a,b,c \in X$, we have
$ (a*b)*c=(a*c)*(b*c). $

A {\it rack} is a set with a binary operation that satisfies
(II) and (III).

A typical example of a quandle is a group $X=G$ with
$n$-fold conjugation
as the quandle operation: $a*b=b^{-n} a b^n$.
Racks and quandles have been studied in, for example,
\cite{Brieskorn,FR,Joyce,K&P,Matveev}.

\begin{figure}
\begin{center}
\mbox{
\epsfxsize=3in
\epsfbox{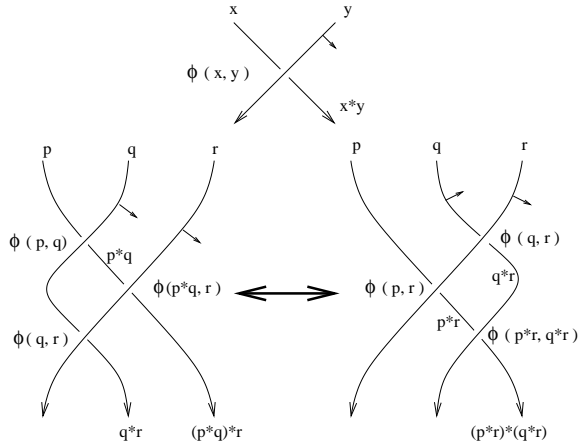}
}
\end{center}
\caption{ Type III move and the quandle identity  }
\label{nu2cocy}
\end{figure}

The axioms for a quandle correspond respectively to the
Reidemeister moves of type I, II, and III
(see
\cite{FR}, \cite{K&P}, for example).
Indeed, knot diagrams were one of the motivations
  to define such an algebraic structure.
In all of our diagrams orientations
and co-orientations
(normal vectors to the given diagram)
are mutually determined
by a right-hand rule. So tangent plus normal
agrees with the counter-clockwise orientation
of the plane that contains the knot diagram.
At a crossing of a classical knot diagram (in which the arcs are
co-oriented), the under-arc is labeled on one
segment by a quandle element, $x$, and along the
other segment by the quandle product
$x*y$ where $y$ is the quandle element
labeling the over-arc. The co-orientation can be used as a
mnemonic for the multiplication; the arc
towards which the normal of the over-arc points receives the product.
See
the top of
Fig.~\ref{nu2cocy}.
The choice of rack multiplication corresponds
to the Wirtinger relation in the fundamental group: $x*y=y^{-1} x y$.
The labels in the figure  involving   $\phi$,
which are assigned to the crossings,  will be used later.

Oriented knotted surface diagrams are  co-oriented by the same rule
(the orientation of the tangent plus the normal vector defines
the given orientation of
3-space);
conversely,
co-orientations of  knotted surface diagrams determine orientations.
The co-orientation is denoted on the complement of
the branch point set by means of a short normal arrow.
In some of the illustrations
only  orientations are indicated, in some only co-orientations are indicated,
  and in some neither are indicated.

A function $f: X \rightarrow  Y$ between quandles
or racks  is a {\it homomorphism}
if $f(a \ast b) = f(a) * f(b)$
for any $a, b \in X$.

\section{Cohomology  of Quandles and Knot Diagrams}
\label{coho}

We define the homology and cohomology theory for racks and quandles.

  Let $C_n^{\rm R}(X)$ be the free
abelian group generated by
$n$-tuples $(x_1, \dots, x_n)$ of elements of a quandle $X$. Define a
homomorphism
$\partial_{n}: C_{n}^{\rm R}(X) \to C_{n-1}^{\rm R}(X)$ by \begin{eqnarray}
\lefteqn{
\partial_{n}(x_1, x_2, \dots, x_n) } \nonumber \\ && =
\sum_{i=2}^{n} (-1)^{i}\left[ (x_1, x_2, \dots, x_{i-1}, x_{i+1},\dots, 
x_n) \right.
\nonumber \\
&&
- \left. (x_1 \ast x_i, x_2 \ast x_i, \dots, x_{i-1}\ast x_i, x_{i+1}, 
\dots, x_n) \right]
\end{eqnarray}
for $n \geq 2$
and $\partial_n=0$ for
$n \leq 1$.
  Then
$C_\ast^{\rm R}(X)
= \{C_n^{\rm R}(X), \partial_n \}$ is a chain complex.

Let $C_n^{\rm D}(X)$ be the subset of $C_n^{\rm R}(X)$ generated
by $n$-tuples $(x_1, \dots, x_n)$
with $x_{i}=x_{i+1}$ for some $i \in \{1, \dots,n-1\}$ if $n \geq 2$;
otherwise let $C_n^{\rm D}(X)=0$. If $X$ is a quandle, then
$\partial_n(C_n^{\rm D}(X)) \subset C_{n-1}^{\rm D}(X)$ and
$C_\ast^{\rm D}(X) = \{ C_n^{\rm D}(X), \partial_n \}$ is a sub-complex of
$C_\ast^{\rm
R}(X)$. Put $C_n^{\rm Q}(X) = C_n^{\rm R}(X)/ C_n^{\rm D}(X)$ and
$C_\ast^{\rm Q}(X) = \{ C_n^{\rm Q}(X), \partial'_n \}$,
where $\partial'_n$ is the induced homomorphism.
Henceforth, all boundary maps will be denoted by $\partial_n$.

For an abelian group $A$, define the chain and cochain complexes
\begin{eqnarray}
C_\ast^{\rm W}(X;A) = C_\ast^{\rm W}(X) \otimes A, \quad && \partial =
\partial \otimes {\rm id}; \\ C^\ast_{\rm W}(X;A) = {\rm Hom}(C_\ast^{\rm
W}(X), A), \quad
&& \delta= {\rm Hom}(\partial, {\rm id})
\end{eqnarray}
in the usual way, where ${\rm W}$
  $={\rm D}$, ${\rm R}$, ${\rm Q}$.

\begin{sect}{\bf Definition.\/} {\rm
The $n$\/th {\it rack homology group\/} and the $n$\/th {\it rack
cohomology group\/} \cite{FRS1}
of a rack/quandle $X$ with coefficient group $A$ are \begin{eqnarray}
H_n^{\rm R}(X;A) = H_{n}(C_\ast^{\rm R}(X;A)), \quad
H^n_{\rm R}(X;A) = H^{n}(C^\ast_{\rm R}(X;A)). \end{eqnarray}
The $n$\/th {\it degeneration homology group\/} and the $n$\/th
  {\it degeneration cohomology group\/}
  of a quandle $X$ with coefficient group $A$ are
\begin{eqnarray}
H_n^{\rm D}(X;A) = H_{n}(C_\ast^{\rm D}(X;A)), \quad
H^n_{\rm D}(X;A) = H^{n}(C^\ast_{\rm D}(X;A)). \end{eqnarray}
The $n$\/th {\it quandle homology group\/}  and the $n$\/th
{\it quandle cohomology group\/ }
of a quandle $X$ with coefficient group $A$ are
\begin{eqnarray}
H_n^{\rm Q}(X;A) = H_{n}(C_\ast^{\rm Q}(X;A)), \quad
H^n_{\rm Q}(X;A) = H^{n}(C^\ast_{\rm Q}(X;A)). \end{eqnarray}

\begin{sloppypar}
The cycle and boundary groups
(resp. cocycle and coboundary groups)
are denoted by $Z_n^{\rm W}(X;A)$ and $B_n^{\rm W}(X;A)$
(resp.  $Z^n_{\rm W}(X;A)$ and $B^n_{\rm W}(X;A)$),
  so that
$$H_n^{\rm W}(X;A) = Z_n^{\rm W}(X;A)/ B_n^{\rm W}(X;A),
\; H^n_{\rm W}(X;A) = Z^n_{\rm W}(X;A)/ B^n_{\rm W}(X;A)$$
where ${\rm W}$ is one of ${\rm D}$, ${\rm R}$, ${\rm Q}$.
We will omit the coefficient group $A$ if $A = {\bf Z}$ as 
usual.\end{sloppypar}

Here we are almost exclusively interested in quandle homology or cohomology.
So
we drop the superscript/subscript
  ${\rm W}= {\rm Q}$ from the notation,
  unless it is needed.
}\end{sect}

\begin{sect} {\bf  Remark.\/} {\rm
Recall that
$C_n^{\rm D}(X; A )$ is the
subgroup of $C_n^{\rm R}(X; A )$ generated by
$\vec{x}=(x_1, \cdots, x_n) \in C_n^{\rm R}(X;A)$ such that
$x_j=x_{j+1}$ for some $j = 1, \cdots, n-1$.
Let $P^n(X ; A)=\{ f \in C^n_{\rm R}(X ; A) | f(\vec{x})=0$
for all $\vec{x} \in C_n^{\rm D}(X) \}$;
this set can be identified with $C^n_{\rm Q}(X ; A)$.
(The set $P^3$ is related to
branch points of knotted surface diagrams.)
Then the quandle cohomology group is described as
  $$H^n_{\rm Q}(X;A) =
  (P^n(X ; A) \cap Z_{\rm R}^n(X; A))/ \delta ( P^{n-1}(X; A)) .$$
There is another cohomology group
defined by
$$H^n_{\rm Q}(X;A)' 
=  (P^n(X;A) \cap Z_{\rm R}^n(X;A))/
(P^n(X;A) \cap B_{\rm R}^n(X;A)).$$
This cohomology group makes sense even for $X$ a rack.
Studies of
this
cohomology
group, 
in relation to the
cohomology
group 
$H^n_{\rm Q}(X;A)$ 
and branch points, are expected.
} \end{sect}

\begin{figure}
\begin{center}
\mbox{
\epsfxsize=3in
\epsfbox{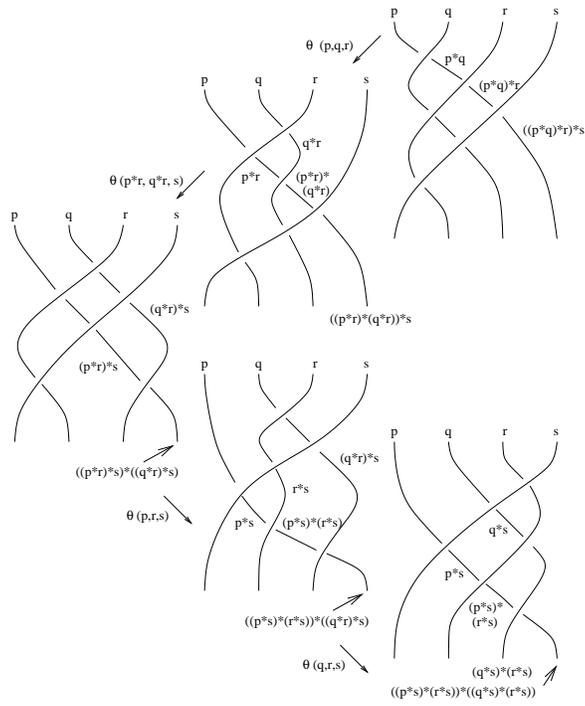}
}
\end{center}
\caption{The tetrahedral move and a cocycle relation, LHS  }
\label{tetraL}
\end{figure}

\begin{figure}
\begin{center}
\mbox{
\epsfxsize=3in
\epsfbox{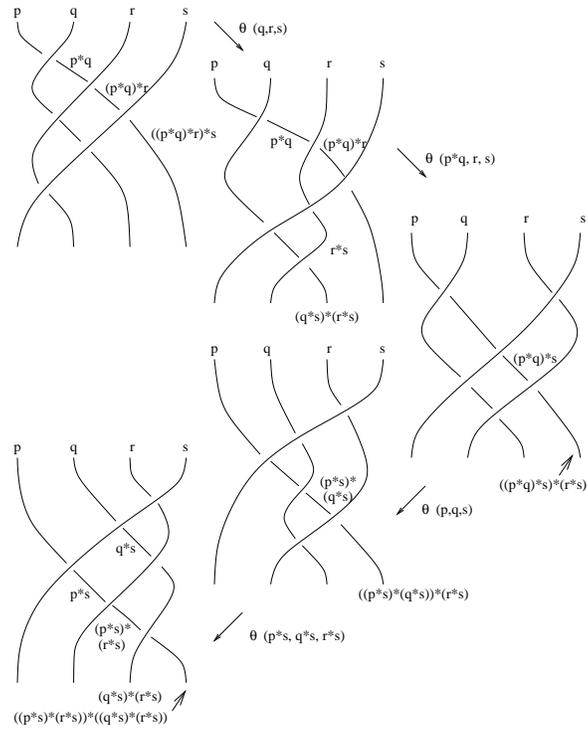}
}
\end{center}
\caption{The tetrahedral move and a cocycle relation, RHS   }
\label{tetraR}
\end{figure}

\begin{sect}{\bf Examples.\/ }{\rm The cocycle conditions
are related to moves on knots and higher dimensional knots as
indicated in Figs.~\ref{nu2cocy},
\ref{tetraL} and \ref{tetraR}.
A 2-cocycle $\phi$ satisfies the relation:
$$\phi(p,r) + \phi(p*r,q*r) = \phi(p,q) +\phi(p*q,r).$$
And a 3-cocycle $\theta$ satisfies the relation:
\begin{eqnarray*}
\lefteqn{ \theta (p,q,r) + \theta (p*r, q*r, s) +
\theta (p,r,s) } \\
& = & \theta (p*q, r, s) + \theta(p,q,s) + \theta (p*s, q*s, r*s) .
\end{eqnarray*}
In subsequent sections, such cocycles will be
assigned to crossings of classical diagrams or triple points
of knotted surface diagrams, respectively.
Figure~\ref{nu2cocy}
shows that
the sum of cocycles evaluated on quandle elements around the
crossings of a diagram remains invariant under
a Reidemeister type III move.
The corresponding
move for knotted surfaces (right-bottom
of Fig.\ref{rose}), called the {\it tetrahedral move}, with
choices of a height function and crossing information, is depicted in
Figs.~\ref{tetraL} and \ref{tetraR}.
Although such  figures involve four straight lines as cross sections
of four planes in space, in Figs.~\ref{tetraL} and \ref{tetraR}
we depicted curved lines instead, to make the figures look nicer.
A $3$-cocycle is assigned to each
type III move in the figures;
these moves correspond
to triple points of a knotted surface diagram.
Thus the sum of $3$-cocycles
(evaluated on the quandle elements near the triple point)
  remains invariant under this move.
Hence the
cocycles can be used to define knot invariants.
We turn now to a rigorous definition of such invariants.
} \end{sect}

\section{Cocycle Invariants of Classical Knots}
\label{classical}

\begin{sect} {\bf Definition.\/} {\rm
A  {\it coloring}
on an oriented  classical knot diagram is a
function ${\cal C} : R \rightarrow X$, where $X$ is a fixed
quandle
and $R$ is the set of over-arcs in the diagram,
satisfying the  condition
depicted
in the top
of Fig.~\ref{nu2cocy}.
In the figure, a
crossing with
over-arc, $r$, has color ${\cal C}(r)= y \in X$.
The under-arcs are called $r_1$ and $r_2$ from top to bottom;
the normal of the over-arc $r$ points from $r_1$ to $r_2$.
Then it is required that
${\cal C}(r_1)= x$ and ${\cal C}(r_2)=x*y$.

Note that locally the colors do not depend on the
orientation of the under-arc.
The quandle element ${\cal C}(r)$ assigned to an arc $r$ by a coloring
  ${\cal C}$ is called a {\it color} of the arc.
This definition of colorings on knot diagrams has been known, see
\cite{FR,FoxTrip} for example.

{\sc Henceforth, all the quandles that are used to color diagrams will be 
finite.}
} \end{sect}

\begin{figure}
\begin{center}
\mbox{
\epsfxsize=4in
\epsfbox{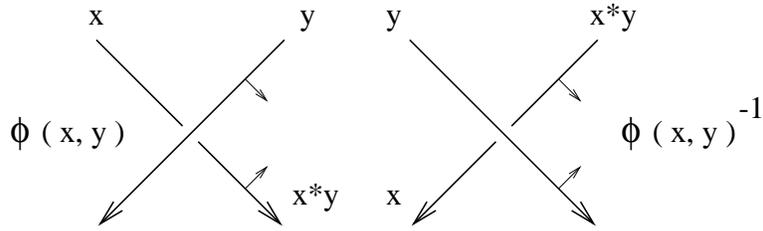}
}
\end{center}
\caption{Weights for positive and negative crossings  }
\label{twocrossings}
\end{figure}

At a crossing,
if the pair of the co-orientation
  of the
over-arc and  that of the under-arc
matches the (right-hand) orientation of the plane, then the
crossing is called {\it positive}; otherwise it is {\it negative}.
In Fig.~\ref{twocrossings},
the two possible oriented and co-oriented
crossings are depicted. The left is a positive crossing,
and the right is negative.

In what follows in this section, we suppose that
a finite quandle $X$ which is used for colorings and
an abelian coefficient group $A$ are fixed.

\begin{sect} {\bf Definition.\/} {\rm
Let $\phi \in Z^2_{\rm Q}(X;A)$ be a $2$-cocycle. 
A ({\it Boltzmann}) {\it weight}, $B(\tau, {\cal C})$,
(associated with $\phi$) 
at a  crossing $\tau$ is defined as follows.
Let  ${\cal  C}$
denote a coloring.
Let $r$ be the over-arc at $\tau$, and $r_1$, $r_2$ be
under-arcs such that the
normal
to $r$ points from $r_1$ to $r_2$.
Let $x={\cal C}(r_1)$ and $y={\cal C}(r)$.
Then define $B(\tau, {\cal C})= \phi(x,y)^{\epsilon (\tau)}$,
where
$\epsilon (\tau)= 1$ or $-1$, if  the sign of $\tau$
is positive or negative, respectively.
} \end{sect}

\begin{sect} {\bf Definition.\/} {\rm
Let $\phi \in Z^2_{\rm Q}(X;A)$ be a $2$-cocycle. 
The {\it partition function}, or a {\it state-sum},
(associated with $\phi$) of a knot diagram  
is the expression
$$
\sum_{{\cal C}}  \prod_{\tau}  B( \tau, {\cal C}).
$$
The product is taken over all crossings of the given diagram,
and the sum is taken over all possible colorings.
(The value of  $  B( \tau, {\cal C})$ is  in the coefficient group $A$ written
multiplicatively). 
The formal sum is taken over all
colorings, and hence
the values of the state-sum
are  in  the group ring
${\bf Z}[A]$. 
} \end{sect}

\begin{figure}
\begin{center}
\mbox{
\epsfxsize=2.5in
\epsfbox{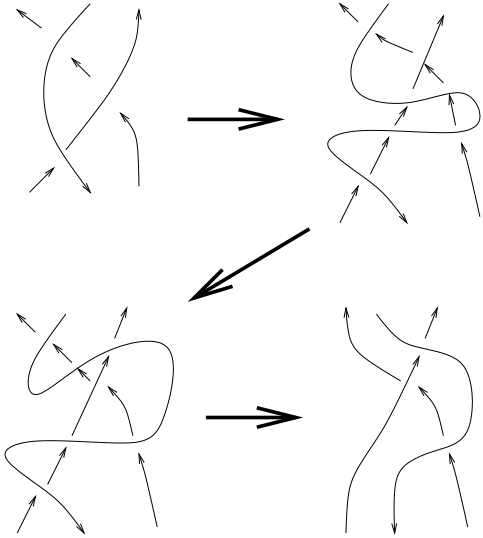}
}
\end{center}
\caption{
A type III move with different crossings  }
\label{kandptrick}
\end{figure}

\begin{sect} {\bf Theorem.\/}
Let $\phi \in Z^2_{\rm Q}(X;A)$ be a $2$-cocycle. 
The partition function
(associated with $\phi$) of a knot diagram 
is invariant under Reidemeister moves,
so that it defines an invariant of knots and links.
Thus it will be  denoted by $\Phi (K)$
(or $\Phi_{\phi}(K)$ to specify the $2$-cocycle $\phi$ used).
\end{sect}
{\it Proof.}
There is a one-to-one correspondence
between colorings before and
after each Reidemeister move. Hence we
check that the state-sum remains
unchanged under Reidemeister moves for
each coloring.
For the type I move, the weight assigned
to the crossing is of
the form $\phi(x,x)^{\pm 1}$, which is $1$
by assumption that $\phi$ is
a quandle cocycle. Thus the state-sum is
invariant under type I moves.
Recall from \cite{K&P}, there are two types of
type II moves depending on whether the
arcs are oriented  in the same direction
or different directions.
  In either case,
at the two crossings of  a type II move,
the 2-cocycle weights
are the same, but with
opposing exponents. Therefore the weights cancel
in the state-sum and the partition function is invariant under type II
moves.

The definition of cocycles was formulated so that the
partition function would be invariant under the
type III move depicted in Fig.~\ref{nu2cocy}.
There are other possible type III moves
depending on the signs
of the crossings and
the  orientation of the edges of the central triangle.
In \cite{K&P}
page 81, Kauffman presents the sketch of the argument
which
shows that the
type III move with differently oriented
triangle follows from the type II moves and
one  choice of type III move.
In Fig.~\ref{kandptrick},
we indicate how to change
the sign of one of the crossings via an analogous
technique. This shows
invariance under all type III moves.
See also \cite{Turaev} or \cite{turaevbook}.
  $\Box$

\begin{sect} {\bf Proposition.\/} \label{coblemma3}
Let $\phi, \phi' \in Z^2_{\rm Q}(X;A)$ be $2$-cocycles. 
If $\Phi_{\phi}$ and $\Phi_{\phi '} $ denote the state-sum invariants
defined from cohomologous cocycles  $\phi$ and $\phi'$
(so that $\phi = \phi' \delta \psi$ for some $1$-cochain $\psi$),
then $\Phi_{\phi} =\Phi_{\phi '} $ (so that $\Phi_{\phi} (K)=\Phi_{\phi '}(K)$
  for any link $K$).

In particular,
the state-sum is equal to the number of
colorings of
a given
knot diagram
if the $2$-cocycle used for the Boltzmann weight is a coboundary.
\end{sect}
{\it Proof.}
We prove the second half,
as the first half follows
from a similar argument.
Suppose that
$\phi (x, y)= \psi(x) \psi(x*y)^{-1}$, so that the cocycle
is a coboundary.  Pick a coloring of the diagram. We can think of the
weight
as a weight
  of the ends of each under-arc where the ``bottom''
  end
of a positive crossing
receives the weight
$\psi(x)$ while the top end of the under-arc
  receives a weight of
$\psi(x*y)^{-1}$.
The negative-crossing  case is  similar.
The under-arc has only one color associated to it,
  so for a given arc, the weights at its two ends cancel.
  A given term in the state-sum then contributes a term of $1$
  to the state-sum. Thus when we sum over all colorings, we end up
counting the colorings.
$\Box$

\vspace{5mm}
We say that the state-sum invariant of a knot/link $K$ is {\it trivial} if 
it is an integer.
In this case, the integer is equal to the number of colorings of a diagram 
of $K$
by $X$.

\section{Cocycle Invariants of Knotted Surfaces}
\label{surface}

First we recall the
notion of
knotted surface diagrams.
See \cite{CS:book} for details and examples.
Let $f:F \rightarrow {\bf R}^4$ denote a smooth embedding of a closed
surface $F$ into 4-dimensional space.
Such an embedding $f$, or its image $f(F)$, is called
{\it a knotted surface}.
By deforming the map $f$ slightly by an ambient isotopy of ${\bf R}^4$
if necessary,
we may assume that
$p \circ f$ is a general position map,
where  $p: {\bf R}^4 \rightarrow {\bf R}^3$
denotes the
orthogonal
projection onto an affine subspace
which does not intersect $f(F)$.
Along the double curves, one of the sheets (called
the {\it over-sheet}) lies farther than the other ({\it under-sheet})
with respect to the projection direction.
The {\it under-sheets}
  are coherently broken in the projection,
and such broken surfaces are called {\it knotted surface diagrams}.

When the surface is oriented, we take normal vectors $\vec{n}$
to the projection of the surface such that the triple
$(\vec{v}_1, \vec{v}_2, \vec{n})$ matches the orientaion of 3-space,
where $(\vec{v}_1, \vec{v}_2)$ defines the orientation of the surface.
Such normal vectors are defined on the projection at all points other than
the isolated branch points.

\begin{figure}
\begin{center}
\mbox{
\epsfxsize=4in
\epsfbox{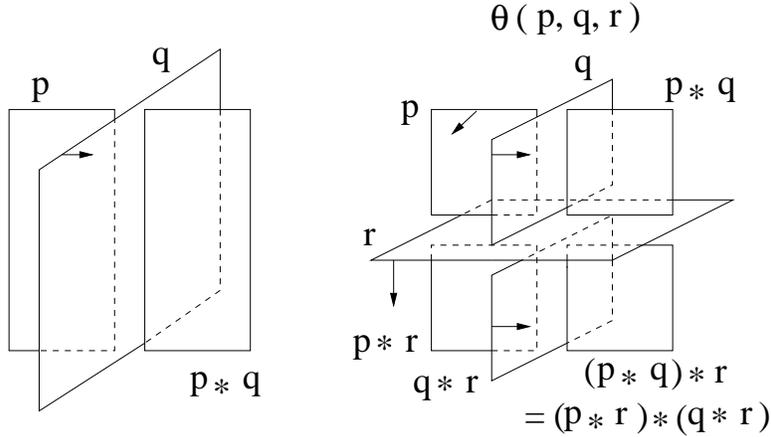}
}
\end{center}
\caption{Colors at a double curve and a triple point   }
\label{triplepoint}
\end{figure}

\vspace{5mm} 
We fix a finite quandle $X$ and an abelian group $A$. 

\begin{sect} {\bf Definition.\/} \label{4dcolor} {\rm
A {\it coloring} on an oriented  (broken) knotted surface diagram is a
function ${\cal C} : R \rightarrow X$, where
$R$ is the set of regions in the broken surface diagram,
satisfying the following condition at the double point set.

At a double point curve, two coordinate planes intersect locally.
One is the
over-sheet $r$, the other is the under-sheet, and the under-sheet is
broken into two components, say $r_1$ and $r_2$.
  A normal of the over-sheet $r$ points to
one of the components, say $r_2$.
If ${\cal C}(r_1) = x \in X$, ${\cal C} (r) = y$, then we require that
${\cal C} (r_2) = x*y$.
The quandle elements ${\cal C} (r)$ assigned to an arc $r$ by a coloring
is called a {\it color} of $r$.
See Fig.~\ref{triplepoint} left.
} \end{sect}

\begin{sect} {\bf Lemma.\/}
The above condition is compatible at each triple point.
\end{sect}
{\it Proof.\/}
The meaning of this lemma is as follows.
There are 6 double curves near a triple point, giving
6 conditions on colors assigned.
  It can be checked in a straightforward manner
that these conditions do not contradict each other.
In particular, there is one of the 4 pieces of the lower
  sheet that receives color
$(a*b)*c$ or $(a*c)*(b*c)$ depending on what path was followed
  to compute the color.
Since these values agree in the quandle, there is no contradiction.
Figure~\ref{triplepoint} illustrates the situation.
$\Box$

\begin{sect} {\bf Definition.\/}
\label{triplepointsign}
{\rm
Note that when three sheets form a triple point, they have relative
positions {\it top, middle, bottom}
with respect to the projection
direction of $p: {\bf R}^4 \rightarrow {\bf R}^3$.
The {\it sign of a triple point}
is positive
if the normals of top, middle, bottom sheets in this order
match the orientation of the $3$-space. Otherwise
the sign is negative.
We use the right-hand rule convention for the orientation of $3$-space.
This definition is found, for example, in \cite{CS:book}.
} \end{sect}

\begin{sect} {\bf Definition.\/} {\rm
Fix a $3$-cocycle $\theta \in Z^3_{\rm Q}(X;A)$.  
A ({\it Boltzmann}) {\it weight}  at a triple point, $\tau$,
is defined as follows.
Let $R$ be the octant from which all normal vectors of the
three sheets point outwards; let
a  coloring ${\cal C}$ be given.
Let $p$, $q$, $r$ be colors of the
bottom, middle, and top
sheets respectively, that bound the region $R$.
Let  $\epsilon (\tau) =1 $ or $-1$ if $\tau$ is positive or negative,
  respectively.
Then the Boltzman weight $ B( \tau, {\cal C})$
at $\tau$ with respect to ${\cal C}$
is defined to be $\theta (p,q,r) ^{ \epsilon (\tau) }$
where $p$, $q$, $r$ are colors described above.
Figure~\ref{triplepoint} illustrates the situation.
} \end{sect}

\begin{sect} {\bf Definition.\/} {\rm
Let $\theta \in Z^3_{\rm Q}(X;A)$ be a $3$-cocycle. 
The {\it partition function}, or a {\it state-sum},
(associated with $\theta$) of a knotted surface diagram  
is the expression
$$ \sum_{\cal C}  \prod_{\tau}  B( \tau, {\cal C} ). $$
The product is taken over all triple points of the diagram,
and the sum is taken over all possible colorings. 
As in the classical case,
$A$ is written multiplicatively and the state-sum is
an element of the group ring ${\bf Z}[A]$. 
} \end{sect}

\begin{figure}
\begin{center}
\mbox{
\epsfxsize=5in
\epsfbox{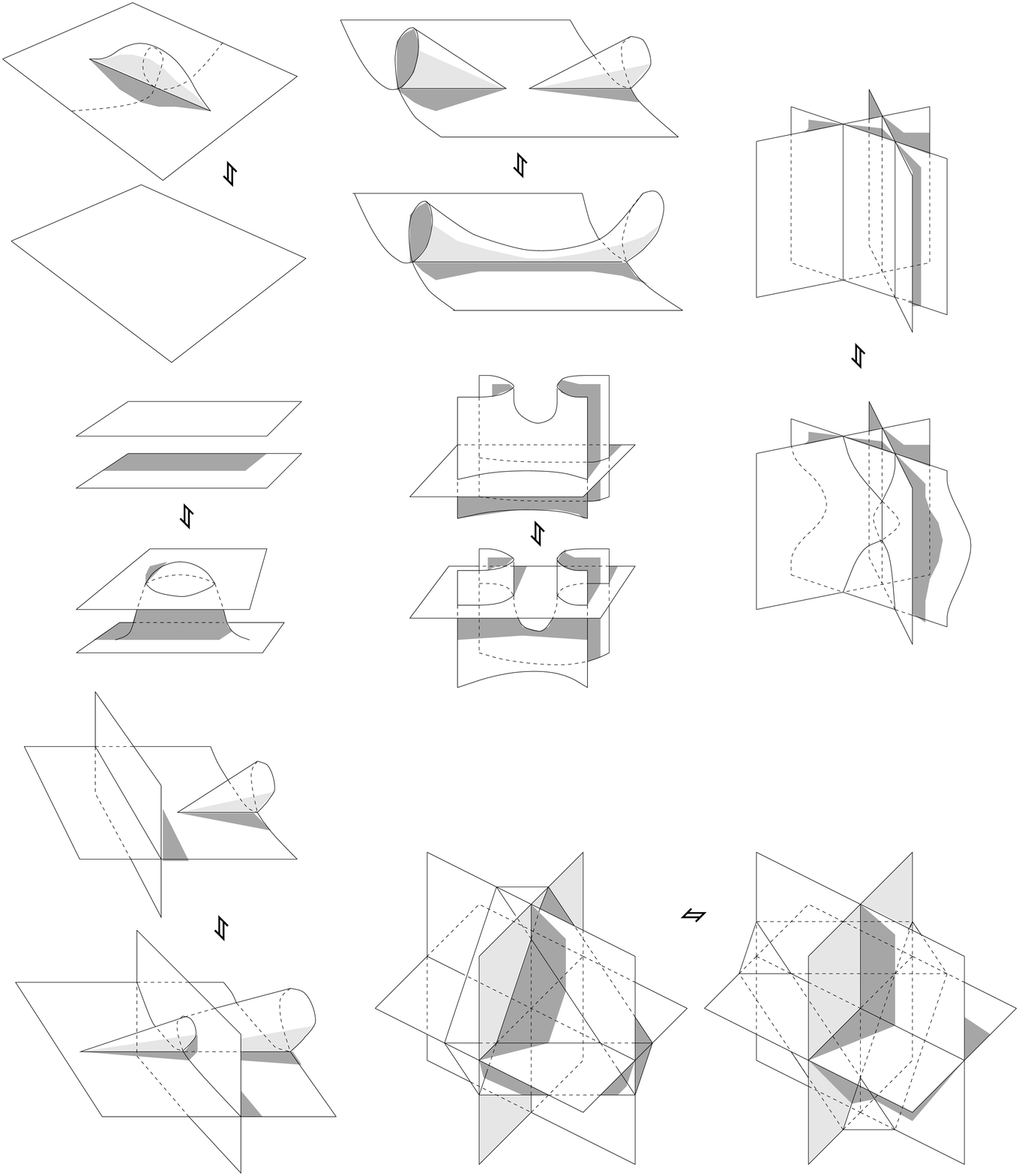}
}
\end{center}
\caption{Roseman moves for knotted surface diagrams  }
\label{rose}
\end{figure}

\begin{figure}
\begin{center}
\mbox{
\epsfxsize=3in
\epsfbox{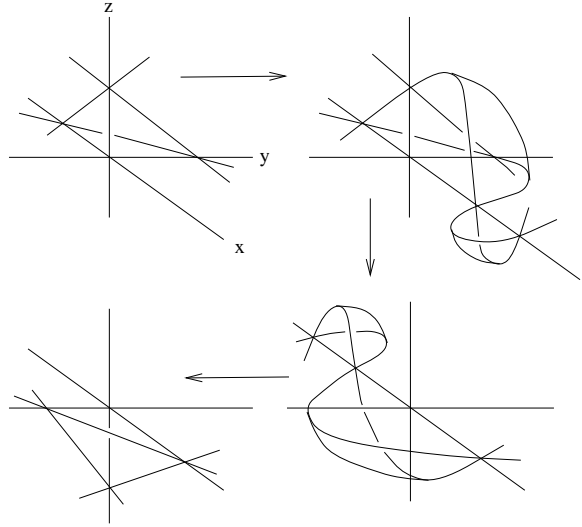}
}
\end{center}
\caption{Turaev's technique generalized to four dimensions  }
\label{turaev4D}
\end{figure}

\begin{sect} {\bf Theorem.\/}
We fix a $3$-cocycle $\theta \in Z^3_{\rm Q}(X;A)$.
The partition function does not depend on the choice of
knotted surface diagram.
Thus it is an invariant of knotted surfaces $F$,
and denoted by $\Phi (F)$
(or $\Phi_{\theta} (F) $ to specify the $3$-cocycle $\theta$ used).
\end{sect}
{\it Proof.\/}
Roseman
provided analogues of the
Reidemeister moves
as moves to
knotted surface
diagrams
and these
analogous
moves (called Roseman moves)
are depicted in Fig.~\ref{rose} \cite{CS:book,Rose}.
Note that in this figure projections are depicted, instead
of broken surface diagrams. There are moves for all possible crossing
information for the sheets involved in each move.
Thus two knotted surface diagrams represent isotopic
knotted surfaces
if and only if the diagrams are related by a finite sequence
of moves,
called {\it Roseman moves,}
taken from this list.
There is a one-to-one correspondence between colorings
before and after each Roseman move, so we check
that the state-sum is invariant under each  Roseman move.
  However, the state-sum depends only
on triple points, so we need only
  consider those moves that involve triple points. These moves are:
(1) the
creation or cancelation of a pair
of oppositely signed triple points
(indicated in the last column of the
second row of the illustration);
(2) moving a branch point through a sheet
(third row, first column);  and  (3) the
tetrahedral move
(on the bottom right)
that motivated
the definition of the cocycles.

In the first case, the pair of triple
points have opposite signs, so for a given
coloring, the two contributing factors
of the state-sum cancel.
In the second case, the branch point
occurs on either the bottom/middle sheet
or on the top/middle sheet,
and these
sheets have the same color. Since the
weighting of the proximate triple point
is a quandle cocycle
(so $\theta(x,x,y)=\theta(x, y, y) =1$),
this factor does not contribute to the state-sum.

In the third case, there are several possible
tetrahedral moves to consider
that depend on (a) the local orientation
of the sheets around the tetrahedron, and
(b) the signs of the triple points that are
the vertices of the tetrahedron.
The definition of the cocycles and the
illustrations Fig.~\ref{tetraL} and Fig.~\ref{tetraR}
indicate that the state-sum is invariant
under  one  of these possible choices.
We will move a given tetrahedral move so
that the planes involved coincide with planes in this
standard position, but have possibly
differing crossings or orientations.
Then we generalize Turaev's technique to
dimension 4 to show that the given move follows from
the fixed move and invariance under adding
or subtracting a cancelling pair of triple points.

Let $T_0$ be the four planes together with
the choice of orientations
depicted in Fig.~\ref{tetraL} and \ref{tetraR},
and let $T$ be a given situation of four planes.
Each sheet has relative  height in $4$-dimensions.
Call them $1$st, $2$nd, $3$rd, $4$th sheet
from bottom to top,
  respectively, so that the
$1$st sheet is the bottom sheet and the $4$th is top.
In other words, the $4$th sheet is unbroken, and the $1$st sheet is broken
into seven pieces in the broken surface diagram.
Suppose that $T_0$ has $xy$, $xz$, $yz$ planes
  as $4$th, $3$rd, and the $2$nd sheets respectively,
and $x+y+z=1$ as the $1$st plane.
We can isotope $T$ to $T_0$ in such a way that
the $4$th sheets match with orientations.
(In other words, isotope the $4$th sheet of $T$ to that of $T_0$
so that the co-orientation normals also match.)
By further isotopy match the $3$rd sheets together with their normals.
Then match the $2$nd sheets. However, the normals may not match here
(if the sign of the triple point among the sheets $2$, $3$, $4$ are opposite).
By isotopy, the $1$st sheet of $T$ is one of the planes
$\pm x \pm y \pm z =1$. However, there are two cases (before/after)
of positions for a given tetrahedral move, so that we may assume
that the $1$st plane is one of four : $\pm x \pm y + z = 1$.
Thus we have four possibilities of orientation choices, those for $2$nd and
the $1$st sheets, and four possibilities for the position of the bottom
($1$st) sheet.

Consider the case where the $1$st sheet has the opposite orientation
of the fixed situation in $T_0$.
Then comparing to $T_0$, the signs of the triple points
involving the $1$st sheet
reverse, and those for the triple point not having the sheet $1$ remain
the same. The former are the triple points among the sheets
$(123)$, $(134)$, $(124)$, and the latter is $(234)$. However recall that
the cocycle assigned to the triple point $(234)$ is $\theta (q,r,s)$
and appears in both sides of the move, and the rest of terms are
inverses of the $2$-cocycle condition.
Therefore the state-sum remains invariant in this case.
Hence the state-sum does not depend on the choice of the orientation
of the $1$st sheet.

Suppose the $1$st sheet is $-x+y+z=1$ as depicted in
the top left of Fig.~\ref{turaev4D}. Then
the
figure shows that
this case follows from the tetrahedral move with the $1$st sheet
$x+y+z=1$,
together with creation/cancelation of a pair of triple points.
  Repeating this process, we conclude that the state-sum is
independent of choice of the position of the $1$st sheet.

It remains to prove that the state-sum is independent of choice of
the orientation of the $2$nd sheet. In other words,
we check the case where the sign of
the triple point among the sheets $4$, $3$, $2$ is negative.
In Fig.~\ref{turaev4D},
regard the $xy$-plane, $xz$-plane, and the plane $-x+y+z=1$
as the sheet $4$, $3$, and $2$ respectively.
Then the tetrahedral move from top right to the bottom right
of the figure is the one with the opposite sign for the
the sheets $4$, $3$, and $2$, comparing to the tetrahedral move
from top left to bottom left.
Therefore this case also does not change the state-sum.
This   completes the proof. $\Box$

\begin{sect} {\bf Proposition.\/} \label{coblemma4}
If $\Phi_{\theta}$ and $\Phi_{\theta '}$ denote the state-sum invariants
defined from cohomologous cocycles  $\theta$ and $\theta '$
(so that $\theta = \theta ' \delta \psi$ for some $2$-cochain $\psi$),
then $\Phi_{\theta} =\Phi_{\theta '} $
(so that $\Phi_{\theta} (K)=\Phi_{\theta '}(K)$ for any knotted surface $K$).

In particular, the state-sum is equal to the number of colorings
of a given knotted surface diagram if the $3$-cocycle used for
the Boltzmann weight is a coboundary.

\end{sect}
{\it Proof.}
We prove the second statement; the first follows
  a similar argument.
Any orientable knotted surface can be
isotoped in 4-dimensional space so that the projection
has no branch points \cite{giller}
(see also
\cite{CS:cancel}).
Thus we assume that the given projection does not have branch points.
Then the double point set of the projection is a graph with $6$-valent
vertices (that are triple points of the projection) possibly
with loops with no vertices.
Now the graph can be directed as follows \cite{CS:book}.
Give an orientation, specified by a vector $\vec{v}$,
  of an edge $e$ in such a way that $\vec{v}$
together with  the normals $\vec{n_1}$,
$\vec{n}_2$ of the top and bottom sheets of the projection
(the triple $(\vec{v}, \vec{n_1}, \vec{n_2} )$ )
  matches the orientation of the three-space.
There are three edges
going into
the triple point and
three
coming out.

Suppose the $\theta$ is a coboundary.
Then it is written as
$$ \theta(p,q,r) = \phi(p,r)^{-1} \phi(p, q)
\phi(p*q, r) \phi (p*r, q*r)^{-1} \phi(q,r )^{-1} \phi (q, r) $$
where the canceling pair is added at the end of the RHS.
The expression on the RHS consists of $2$-cochains
$\phi( x, y)$ where the pairs $(x, y)$ ranges over
all pairs that appear at double curves near a triple point,
when the triple point is colored in such a way that $\theta (p,q,r)^{\pm 1}$
is its weight, where $\pm 1$
is the sign of the triple point.
In other words the cochain $\phi(x,y)$ corresponds to an edge with colors
$x$, $y$, $x*y$.
Furthermore, if an edge is oriented out of the triple point,
then the cochain has
negative
exponent $\phi (x, y)$,
and if the edge is oriented into the triple point,
then the cochain has
positive
  exponent. This is true
for either sign ($\pm 1$) of the triple point.

The weight $\theta (p,q,r)^{\pm 1}$ may be considered
to be the  values
$\phi (x, y)^{\pm 1}$  assigned to
  the end points of six edges at the triple point.
Note that if the edge comes out of the triple point, then
the end point receives $\phi(x,y)^{-1}$, and if it comes in,
then it receives $\phi(x,y)$.
Now the state-sum for a fixed color  is the product
of all these weights assigned to end points of oriented edges.
However, for every edge, the initial end point receives
$\phi (x,y)^{-1}$ and the terminal point receives
$\phi (x,y)$, so that the product of weights cancel out.
Thus the product is $1$ for any color, hence the state-sum
is the number of colors of the diagram.
$\Box$

\begin{sect} {\bf Remark.\/} {\rm
The fundamental quandle is defined (see \cite{Joyce} for example)
for codimension $2$ embeddings, and its presentations
are defined (see \cite{FR} for example) using knot diagrams
in all dimensions, along the line of Wirtinger presentations
of knot groups.
Let $Q(K)$ be the fundamental quandle of a knotted curve or surface $K$,
and $X$ be a finite quandle. Then it is seen using presentations
of $Q(K)$ that there is a one-to-one correspondence between
  quandle homomorphisms $Q(K) \rightarrow X$ and
colorings ${\cal C}: R \rightarrow X$.

} \end{sect}

\section{Computing Quandle Cohomology}
\label{quandleco}

In this section, quandle cohomology groups will be
computed for
some
interesting examples.
The cohomology groups will be computed directly from the definitions.
More advanced techniques, such as exact sequences, would be desirable,
but are not pursued here.
We present some computational details, as
some of the calculations will be
used in later sections
  to find non-trivial invariants.

\begin{sect}{\bf Remark.\/}
{\rm Suppose that the coefficient group $A$ is
a cyclic group
written additively:
${\bf Z}$ or
${\bf Z}_n$.
Define a {\it characteristic function}
$$\chi_{x}(y) = \left\{ \begin{array}{lr} 1 & {\mbox{\rm if}} \  \ x=y
\\
                                                0   & {\mbox{\rm if}} \ \  x\ne
y \end{array} \right.$$
from the free abelian group generated by $X^n$ to the group $A.$

The set $\{ \chi_x
: x \in X^n
\}$ of such
functions
spans the group $C^n_{\rm R}(X;A)$ of cochains.
Thus if  $f \in C^n_{\rm R}(X;A)$ is a cochain, then
$$f = \sum_{x \in X^n}
C_x
\chi_x.$$
If  $f \in C^n_{\rm Q}(X;A)$, then $f$ is written as
$$f = \sum_{x \in X^n \setminus S}     C_x \chi_x, $$
where $S =\{ (x_1, \ldots , x_n) : x_j= x_{j+1}
\ {\mbox{\rm for some}} \  j= 1, \ldots , n-1 \}. $
If $\delta f=0$, then $f$ vanishes on expressions of the form
$$\sum_j (-1)^{j+1}(x_0, \ldots , \hat{x}_j, \ldots, x_n)
+
\sum_k (-1)^k(x_0 * x_k, \ldots, x_{k-1}* x_k, x_{k+1}, \ldots , x_n).$$

In computing the cohomology we consider all such expressions as
$(x_0, \ldots x_n)$ ranges over all $(n+1)$-tuples for which
each consecutive pair of elements is distinct.
By evaluating linear
combinations of
characteristic functions on these expressions, we determine
those functions that are cocycles.
Similarly, we compute the coboundary on each of the characteristic
functions in the previous dimension, to determine which linear
combinations of characteristic functions are coboundaries.
Since $A$
is a cyclic group,
the generator will be denoted $1$ (resp. $t$),
  the identity is denoted $0$ (resp. 1), and the characteristic functions
take values
$0$ or $1$ (resp. 1 or $t$) when $A$ is written additively (resp.
multiplicatively).
We turn now to examples.
All cohomology groups are quandle ones $H^n_{\rm Q}$,
unless otherwise stated.
}\end{sect}

\begin{sect} {\bf Definition \cite{FR}.\/} {\rm
A rack is called {\it trivial} if $x*y=x$ for any $x,y$.

The {\it dihedral quandle} $R_n$ of order $n$ is the quandle consisting of
reflections of the regular $n$-gon with the conjugation as operation.
The dihedral group $D_{2n}$  has a presentation
$$ \langle x, y | x^2 = 1 = y^n, xyx=y^{-1} \rangle $$
where $x$ is a reflection and $y$ is a rotation of a regular $n$-gon.
The set of reflections $R_n$ in this presentation
is $\{ a_i = xy^i : i=0, \cdots, n-1 \} $
where we use the subscripts from ${\bf Z}_n$ in the following computations.
The operation is
$$a_i * a_j = a_j^{-1} a_i a_j =x y^{j} x y^i x y^j
= x y^j y^{-i} y^j = a_{2j-i}.$$
Hence $R_n$ is identified with ${\bf Z}_n=\{0, \cdots, n-1 \}$,
with quandle operation $i*j = 2j-i $ (mod $n$).
Compare with
the well known $n$-coloring of
  knot diagrams \cite{FoxTrip}.

Let $S_4$ denote the quandle
with  four elements,
denoted by  $0,1,2,3,$
with the relations
$$\begin{array}{ccccccccl}0 &=& 0 * 0 &= & 1*2 &=& 2* 3 &=& 3*1 \\
                         1  &=& 0*3 &=& 1*1 &=& 2* 0 &=& 3*2 \\
                          2 &=& 0*1 &=& 1*3&=& 2*2 &=& 3*0 \\
                          3 &=& 0*2 &=& 1*0&=&  2*1&=& 3*3. \end{array}$$
This quandle is the set of
clockwise rotations of the faces of a tetrahedron with
conjugation as the operation.
} \end{sect}

\begin{sect} {\bf Definition \cite{FR,K&P}.\/} {\rm
Let $\Lambda = {\bf Z}[T, T^{-1}]$ be the Laurent polynomial ring
over the integers. Then any $\Lambda$-module $M$
  has a quandle structure defined by
$a*b= Ta + (1-T) b$ for $a, b \in M$.
} \end{sect}

For a  Laurent polynomial $h(T)$
whose leading and terminal coefficients are $\pm 1$,
  ${\bf Z}_n[ T, T^{-1} ]  / (h(T)) $
is  a finite quandle. 
We call such quandles {\it (mod $n$)-Alexander quandles.}
Alexander quandles are of interest in Section~\ref{danssect}.

\begin{sect} {\bf Proposition.\/}
\label{s4islx}
We have the following two quandle isomorphisms:
$$ R_4 \cong  {\bf Z}_2[ T, T^{-1} ] / (T^2-1), \quad
{\mbox{\rm and}} \quad
S_4 \cong {\bf Z}_2[T, T^{-1}]/(T^2 +T +1). $$
\end{sect}
{\it Proof.}
The set of elements of either of these Alexander
quandles can be represented as
$\{ 0, 1, T, 1+T \}$.
The following assignment defines an isomorphism $ R_4 \cong  {\bf Z}_2[ T, 
T^{-1} ] / (T^2-1)$ :
  $0 \leftrightarrow 0$,
       $1 \leftrightarrow  1$,
       $2 \leftrightarrow  1+T, $ and
       $3 \leftrightarrow T$.
It happens that the same correspondence also gives
an isomorphism to $S_4$.
$\Box$

\begin{sect} {\bf Lemma.\/}
Any
cochain
on a  trivial quandle
is a cocycle. Only the zero map is a coboundary.
\end{sect}
{\it Proof.} This follows from the definitions. $\Box$

\vspace{5mm}

It is worth remarking here that the trivial quandle is quite effective in
detecting linking. See Section~\ref{danssect}.

\begin{sect} {\bf Lemma.\/}
$H^2(R_3; {\bf Z} ) \cong 0.$
\end{sect}
{\it Proof.}
Let a  $2$-cocycle $f \in Z^2(R_3; {\bf Z} )$ be expressed as
$$f= \sum_{ i,j \in R_3} C_{(i,j)} \chi_{(i,j)}, $$
then
$$ C_{(p,r)} + C_{(p*r,q*r)} - C_{(p,q)} - C_{(p*q,r)} =0
\quad {\rm for~} p,q,r \in R_3$$
and
$$ C_{(p,p)}= 0 \quad {\rm for~} p  \in R_3.$$

The quandle $R_3$ has three elements, $0,1,2$ with quandle operation
$$i*j = 2j -i \ \ \ {\mbox{\rm (mod $3$)}}.$$
Substituting $0,1,2$ for all possibilities for the variables $p,q,r$
into the above expressions, we have $30$
equations on $C_{(i,j)}$, which are simplified as the following:
\begin{eqnarray*}
C_{(0,1)} + C_{(2,1)}  &=& 0 \\
C_{(0,2)} - C_{(2,0)} + C_{(2,1)}  &=& 0 \\
C_{(1,0)} + C_{(2,0)}  &=& 0 \\
C_{(0,2)} + C_{(1,2)}  &=& 0 \\
C_{(i,i)}    &=&  0 \quad {\rm for~} i  \in \{0,1,2\}. \\
\end{eqnarray*}
Therefore,
$$
\begin{array}{lll}
C_{(0,0)}= 0, & C_{(0,1)}= \alpha, & C_{(0,2)}= \beta, \\
C_{(1,0)}= \alpha-\beta, & C_{(1,1)}= 0, & C_{(1,2)}= -\beta, \\
C_{(2,0)}= \beta-\alpha, & C_{(2,1)}= -\alpha, & C_{(2,2)}= 0, \\
\end{array}
$$
where we put $C_{(0,1)}= \alpha$ and $C_{(0,2)}= \beta$.
Then
$$f= \alpha [ \chi_{(0,1)} + \chi_{(1,0)} - \chi_{(2,0)} - \chi_{(2,1)}  ] +
\beta [ \chi_{(0,2)} - \chi_{(1,0)} - \chi_{(1,2)} + \chi_{(2,0)}  ].$$
Since
\begin{eqnarray*}
\delta \chi_1 &=&
- \chi_{(0,2)} + \chi_{(1,0)} + \chi_{(1,2)} - \chi_{(2,0)} ,  \\
\delta \chi_2 &=&
- \chi_{(0,1)} - \chi_{(1,0)} + \chi_{(2,0)} + \chi_{(2,1)},
\end{eqnarray*}
we see that $f$ is a coboundary.
  $\Box$

\begin{sect} {\bf Lemma.\/}
\label{65}
$H^2( R_4; {\bf Z} ) \cong  {\bf Z} \times {\bf Z} .$
\end{sect}
{\it Proof.}
The quandle $R_4$ has four elements,
$0,1,2,$ and $3$; and quandle  operation is
$$i*j=2j -i \ \ \ {\mbox{\rm (mod $4$) }} $$
Let a  $2$-cocycle $f \in Z^2(R_4; {\bf Z} )$ be expressed as
$$f= \sum_{ i,j \in R_4} C_{(i,j)} \chi_{(i,j)}, $$  
then
\begin{eqnarray*}
C_{(0,1)}-C_{(0,3)}+C_{(2,1)}-C_{(2,3)} &=& 0\\
C_{(0,2)}+C_{(2,1)}-C_{(2,3)} &=&0\\
C_{(1,3)}+C_{(3,1)}& =& 0\\
   C_{(2,0)}+C_{(2,1)}-C_{(2,3)} &=& 0\\
-C_{(2,0)}+C_{(2,1)}-C_{(2,3)} &=&0\\
   C_{(1,0)}-C_{(1,2)}+C_{(3,1)} &=&0\\
-C_{(1,3)}+C_{(3,1)} &=& 0\\
C_{(3,0)}+C_{(3,1)}-C_{(3,2)} &=& 0 \\
C_{(i,i)}    &=&  0 \quad {\rm for~} i  \in \{0,1,2,3\}. \\ 
\end{eqnarray*}

Thus
$$C_{(1,3)}=C_{(3,1)},$$
$$C_{(0,2)}=C_{(2,0)},$$
and
$$2C_{(1,3)}=2C_{(0,2)}=0.$$

The relations among the coeficients give that
the group of 2-cocycles is represented as
$$Z^2( R_4; {\bf Z} ) \cong \Hom ({\bf Z}^4 \times ({\bf Z}_2)^2, {\bf Z} )
  \cong {\bf Z}^4 $$
where the generators
are
\begin{eqnarray*}
f_{(0,1)}  & = & \chi_{(0,1)} + \chi_{(0, 3)} \\
f_{(2,1)} & = & \chi_{(2,1)} + \chi_{(2, 3)} \\
f_{(1,0)}  & = & \chi_{(1,0)} + \chi_{(1, 2)} \\
f_{(3,0)} & = & \chi_{(3,0)} + \chi_{(3, 2)} .
\end{eqnarray*}

The coboundaries are computed as follows.
\begin{eqnarray*}
\delta \chi_{0} &=& \chi_{(0, 1)} + \chi_{(0, 3)}
- \chi_{(2, 1)} - \chi_{(2, 3)} \\
& = &   f_{(0, 1)} - f_{(2, 1)}  \\
\delta \chi_{2} &=& - \chi_{(0, 1)} - \chi_{(0, 3)}
  + \chi_{(2, 1)} + \chi_{(2, 3)} \\
  & = & - f_{(0, 1)} + f_{(2, 1)}  \\
\delta \chi_{1} &=& \chi_{(1, 0)} + \chi_{(1, 2)}
- \chi_{(3, 0)} - \chi_{(3, 2)} \\
  & = &  f_{(1, 0)} - f_{(3, 0)}  \\
\delta \chi_{3} &=& - \chi_{(1, 0)} - \chi_{(1, 2)}
  + \chi_{(3, 0)} + \chi_{(3, 2)} \\
  & = & - f_{(1, 0)} + f_{(3, 0)} .
\end{eqnarray*}

Therefore $H^2(R_4; {\bf Z}) \cong  {\bf Z}^2$.
$\Box$

\vspace{5mm}

We have the following calculations that were performed using
Mathematica and Maple.
(See also \cite{cjks2,Greene}
for more on the 2nd homology of the dihedral quandles).

\begin{sect} {\bf Lemma.\/} For the 3-element dihedral quandle we have:
$$H^3(R_3; {\bf Z}_3) \cong {\bf Z}_3,$$
and
$$H^3(R_3; {\bf Z}) \cong 0.$$
\end{sect}
{\it Proof.\/}
We summarize the calculation. For any coefficient group,
any cocycle can be written as
$\sum_{i=1}^5 a_i \eta_i$ where
$3a_1=0$, and

\begin{eqnarray*}
\eta_{1}&=&-\chi_{(0,1,0)} + \chi_{(0,2,0)}+2\chi_{(0,2,1)}+\chi_{(1,0,1)} 
	+\chi_{(1,0,2)}+\chi_{(2,0,2)}+\chi_{(2,1,2)}; \\
\eta_{2}&=&-\chi_{(0,1,0)}  +\chi_{(0,2,1)}-\chi_{(1,0,1)}+\chi_{(1,2,0)}; \\
\eta_{3}&=&\chi_{(0,1,0)}+\chi_{(0,1,2)}-\chi_{(0,2,0)} -\chi_{(0,2,1)}
				-\chi_{(1,0,2)}+\chi_{(1,2,1)};\\
\eta_{4}&=&\chi_{(0,1,0)}+\chi_{(0,1,2)}-\chi_{(0,2,0)}-\chi_{(0,2,1)}
				+\chi_{(2,0,1)}-\chi_{(2,1,2)};\\
\eta_{5}&=&\chi_{(0,1,2)}-\chi_{(0,2,0)}-\chi_{(2,0,2)}+\chi_{(2,1,0)}.
\end{eqnarray*}

For example, if
the coefficient group is  ${\bf Z}$ , then
$Z^3(R_3; {\bf Z}) = {\bf Z}^4$
and is generated by  $\eta_{2},$ $\eta_{3},$ $\eta_{4},$ and $\eta_{5}.$
If the coefficient group is  ${\bf Z}_3$ , then
$Z^3(R_3; {\bf Z}_3) =({\bf Z}_3)^5 $
and is generated by  $\eta_{1},$ $\eta_{2},$ $\eta_{3},$ $\eta_{4},$  and 
$\eta_{5}.$

The following elements generate the group of coboundaries:
\begin{eqnarray*}
\delta \chi_{(0,1)}
&=&
      -\chi_{(0, 1, 0)} - \chi_{(0, 1, 2)} + \chi_{(0, 2, 0)} + \chi_{(0, 2, 
1)} +
           \chi_{(1, 0, 2)} - \chi_{(1, 2, 1)};\\
\delta \chi_{(0,2)}
&=&
   \chi_{(0, 1, 0)} + \chi_{(0, 1, 2)} - \chi_{(0, 2, 0)} - \chi_{(0, 2, 1)} +
           \chi_{(2, 0, 1)} - \chi_{(2, 1, 2)};\\
\delta \chi_{(1,0)}
&=&
     \chi_{(0, 1, 2)} - \chi_{(0, 2, 0)} - \chi_{(1, 0, 1)} - \chi_{(1, 0, 2)} +
           \chi_{(1, 2, 0)} + \chi_{(1, 2, 1)};\\
\delta \chi_{(1,2)}
&=&
     \chi_{(1, 0, 1)} + \chi_{(1, 0, 2)} - \chi_{(1, 2, 0)} - \chi_{(1, 2, 1)} -
           \chi_{(2, 0, 2)} + \chi_{(2, 1, 0)};\\
\delta \chi_{(2,0)}
  &=&
    -\chi_{(0, 1, 0)} + \chi_{(0, 2, 1)} - \chi_{(2, 0, 1)} - \chi_{(2, 0, 2)} +
           \chi_{(2, 1, 0)} + \chi_{(2, 1, 2)};\\
\delta \chi_{(2,1)}
  &=&
     -\chi_{(1, 0, 1)} + \chi_{(1, 2, 0)} + \chi_{(2, 0, 1)} + \chi_{(2, 0, 
2)} -
           \chi_{(2, 1, 0)} - \chi_{(2, 1, 2)}. \end{eqnarray*}

Comparing the cocycles and coboundaries, we have the result. $\Box$

\begin{sect}{\bf Remark.\ }
{\rm In Section~\ref{themain},
  we  use the $3$-cocycle
$$\eta_1= -\chi_{(0,1,0)} + \chi_{(0,2,0)} -
\chi_{(0,2,1)}+\chi_{(1,0,1)}
	+\chi_{(1,0,2)}+\chi_{(2,0,2)}+\chi_{(2,1,2)}  \in Z^3(R_3;{\bf 
Z}_3)$$
to distinguish the 2-twist spun trefoil from its orientation reversed image.}
\end{sect}

\begin{sect} {\bf Remark.\/} {\rm Similar computations give
the following results that are used to compute knot invariants
in a subsequent paper.
$$ H^2 (S_4; A) = \left\{ \begin{array}{ll}
{\bf Z}_2
& {\rm for~} A= {\bf Z}_2 \\ 
   0
& {\rm for~} A = {\bf Z}
\end{array}\right.$$
$$H^3(S_4;A) = \left\{ \begin{array}{ll} {\bf Z}_2
& {\rm for~} A= {\bf Z} \\
                                          0
& {\rm for~} A = {\bf Q} \\
                            ({\bf Z}_2)^3
& {\rm for~} A ={\bf Z}_2 \\
            ({\bf Z}_2)^2 \times  {\bf Z}_4
& {\rm for~} A = {\bf Z}_4
\end{array}\right.$$
} \end{sect}

\section{Group $2$-cocycles and quandle $2$-cocycles}
\label{reltogroup}

In this section we
give quandle $2$-cocycles
using group $2$-cocycles.
Let $G$ be a group and let $A$ be an abelian group
(written multiplicatively) upon which
the group ring ${\bf Z}[G]$ acts trivially.
Then the group cohomology is defined from the following cochain complex.
The abelian group of all maps from the cartesian product of $n$ copies of
$G$ to $A$ is denoted by $C^n(G; A)$.
A coboundary operator
  $\delta : C^n(G; A) \rightarrow C^{n+1}(G;A)$ is defined by
$$ (\delta f)(x_1, \ldots, x_{n+1} )
= f(x_2, \ldots, x_{n+1}) \prod_{i=1}^n
f(x_1, \ldots, x_i x_{i+1}, \ldots, x_{n+1})^{(-1)^i}
  f(x_1, \ldots, x_n)^{(-1)^{n+1}},$$
where $f \in C^n(G;A)$ and $x_1, \ldots, x_{n+1} \in G$.

\begin{figure}
\begin{center}
\mbox{
\epsfxsize=2.5in
\epsfbox{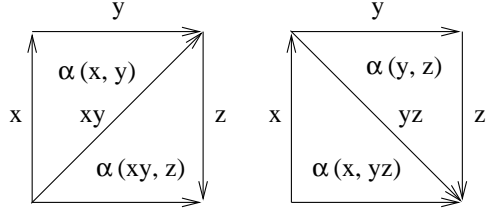}
}
\end{center}
\caption{Group $2$-cocycle condition and triangulations of a square  }
\label{gp2cocy}
\end{figure}

In particular, a function $\alpha: G \times G \rightarrow A$
satisfies {\it the group $2$-cocycle condition}
if the following relation holds:
$$ \alpha(x, y) \alpha(xy, z) = \alpha (x, yz) \alpha(y, z)  .$$
The diagrammatic interpretation of this condition is
depicted in Fig.~\ref{gp2cocy}.
Consider triangulations of planar regions. Suppose the edges are
oriented in such a way that at every triangle, exactly two edges
point to the same orientation (clockwise or counter-clockwise)
and one edge has the opposite direction.
Let $G$ be a finite group, and assign elements of $G$ on the edges,
such that if the two edges of the same directions receive $x$ and $y$ in
this direction, then the other edge receives $xy$.
The value $\alpha(x,y)$ of a 2-cocycle
$\alpha$ is assigned to such a triangle \cite{DW,Fukuma}
(see also \cite{CKS2}).
With this convention, two ways of
triangulating
a square
corresponds to the $2$-cocycle condition as depicted in Fig.~\ref{gp2cocy}.

\begin{figure}
\begin{center}
\mbox{
\epsfxsize=3in
\epsfbox{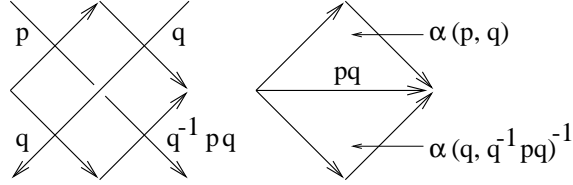}
}
\end{center}
\caption{Defining a quandle cocycle from a group cocycle  }
\label{crossing}
\end{figure}

\begin{sect} {\bf Theorem.\/}
Let $G$ be a group,
considered also as a quandle by conjugation
that we denote by $G_{\rm conj}$.
Let $\alpha \in Z^2(G;A)$ be a group $2$-cocycle.
Define a quandle $2$-cochain
$$\phi (p,q) = \alpha(p,q) \alpha(q, q^{-1} pq)^{-1}.$$
Then $\phi$ is a quandle $2$-cocycle,
$\phi \in  Z^2(G_{\rm conj};A)$.
\end{sect}
{\it Proof.\/}
A
similar argument
to
that in
  \cite{Wakui}
shows that
a group $2$-cocycle $\alpha$ satisfies
\begin{eqnarray*}
\alpha (x,y) & = & \alpha (x^{-1}, xy)^{-1} = \alpha (xy, y^{-1} )^{-1}, \\
  \alpha(x,y) & = &\alpha(z, z^{-1}x) ^{-1} \alpha (z^{-1} x, y)
  \alpha (z, z^{-1}xy ).
\end{eqnarray*}
Using these identities and $2$-cocycle
conditions, one computes
\begin{eqnarray*}
\lefteqn{ \phi (p,q) \phi (p*q, r) \phi (q,r)  } & &  \\
  &=& \alpha (p,q) \alpha (q, q^{-1}pq)^{-1}
\underline{ \alpha ( q^{-1}pq, r) \alpha (r, r^{-1} q^{-1} pqr )^{-1} }
\alpha (q, r) \alpha (r, r^{-1} qr)^{-1} \\
  &=& \alpha (p,q) \underline { \alpha (q, q^{-1}pq)^{-1}
\alpha ( q^{-1}pq, q^{-1} p^{-1} qr) }
\alpha ( q^{-1} p^{-1} qr,  r^{-1} q^{-1} pqr )^{-1}
\underline{ \alpha (q, r) }  \alpha (r, r^{-1} qr)^{-1} \\
  &=& \alpha (p,q)
\underline{ \alpha (pq, q^{-1} p^{-1} qr) }
\alpha ( q^{-1} p^{-1} qr,  r^{-1} q^{-1} pqr )^{-1}
  \alpha (r, r^{-1} qr)^{-1} \\
  &=& \underline{ \alpha (p,q)
\alpha (pr, r^{-1} q) ^{-1} }
\alpha( r^{-1} q,   q^{-1} p^{-1} qr)
\alpha (pr,  r^{-1} p^{-1} qr)
\alpha (q^{-1} p^{-1} qr,  r^{-1} q^{-1} pqr )^{-1}
  \alpha (r, r^{-1} qr)^{-1} \\
  &=& \alpha (p,r) \alpha (r, r^{-1}q )^{-1}
\alpha( r^{-1} q,   q^{-1} p^{-1} qr)
\underline{ \alpha (pr,  r^{-1} p^{-1} qr) }
\alpha (q^{-1} p^{-1} qr,  r^{-1} q^{-1} pqr )^{-1}
\underline{  \alpha (r, r^{-1} qr)^{-1} } \\
  &=& \alpha (p,r) \alpha (r, r^{-1}q )^{-1}
\underline{ \alpha( r^{-1} q,   q^{-1} p^{-1} qr) }
\alpha( r, r^{-1} pr )^{-1} \alpha  ( r^{-1} pr,  r^{-1} p^{-1} qr)
\underline{ \alpha (q^{-1} p^{-1} qr,  r^{-1} q^{-1} pqr )^{-1} } \\
  &=& \alpha (p,r) \underline { \alpha (r, r^{-1}q )^{-1}
  \alpha (r^{-1} q, r)}
  \alpha ( r^{-1} p^{-1} qr,  r^{-1} q^{-1} pqr)^{-1}
\alpha( r, r^{-1} pr )^{-1} \alpha  ( r^{-1} pr,  r^{-1} p^{-1} qr) \\
  &=& \alpha (p,r)
\alpha (q, r) \alpha (r, r^{-1} q r )^{-1}
  \underline{  \alpha ( r^{-1} p^{-1} qr,  r^{-1} q^{-1} pqr)^{-1} }
\alpha( r, r^{-1} pr )^{-1} \underline{ \alpha( r^{-1} pr,  r^{-1} p^{-1} 
qr)}  \\
  &=&\alpha (p,r)
\alpha (q, r) \alpha (r, r^{-1} q r )^{-1}
\alpha ( r^{-1} pr,  r^{-1} q r) \alpha(  r^{-1} q r,  r^{-1} p^{-1} qr)^{-1}
\alpha( r, r^{-1} pr )^{-1} \\
  &=& \phi(q, r) \phi (p,r)
\phi(  p*r, q*r ). \quad \Box
\end{eqnarray*}

The
above
computation
is  easily
carried out using diagrams.
At a crossing, a square is assigned as in Fig.~\ref{crossing} left.
Then the square is triangulated, and group cocycles are assigned
as in Fig.~\ref{gp2cocy}.
The Reidemeister type III move, then, is interpreted as changes
of triangulations
of squares, giving the above computations.

\section{Computations of Cocycle Invariants of  Classical Knots and Links}
\label{danssect}

Suppose a link, $L$, is colored by the trivial $n$-element
quandle, $ T_{n}$,
whose elements we represent by integers $T_n= \{1,\ldots, n  \}$.
  Since $ a * b= a$ for all $ a$ and $b \in T_n$,
each
component of a link $L$
is monochromatically colored.
More precisely, if $a_i$, $i=1, \ldots, m$ are arcs of  a
component $K$ of $L$, the color ${\cal C}(a_i)$ takes the same
value in $T_n$ for $i=1, \ldots, m$.
We assume that the coefficient group $A$ is cyclic
and generated by $t$ (which is infinite cyclic for a while).  
The coboundary homomorphism $\delta$ is trivial for $T_n$ for any $n$,
and  in particular, any function $\phi$ is a cocycle in $ T_{n}$.
Consider the
characteristic functions
(that we write multiplicatively for this section):
$$\chi_{(x,y)}(a,b)= \left\{ \begin{array}{lr} t & {\mbox{\rm if \
}} (a,b) = (x,y), \\ 1  & {\mbox{\rm otherwise.}}
\end{array} \right. $$
For an
  $n$-component link $L=K_1 \cup \ldots \cup K_n$,
let $\lk (K_{i},K_{j})$ denote the linking number of the pair  $(K_{i},K_{j})$
of components, and let  $\lk(L)=\sum_{i<j} \lk(K_i, K_j)$
denote  the {\it total linking number},
where the sum ranges over all pairs with $i<j$, $i,j = 1, \ldots, n$.
Define also the {\it linking number} $\lk (A,B)$
for any disjoint pair of subsets $A, B \subset \{ K_1, \ldots, K_n \}$
  by
$\lk (A, B)= \sum_{ K_u \in A, K_v \in B} \lk(K_u, K_v)$,
where $\lk (A, \emptyset)=0=\lk (\emptyset, B)$.
Recall that the linking number of a $2$-component classical
link $L=K_1 \cup K_2$ can be computed
  by counting the crossing number with signs
($\pm 1$) where the component $K_1$ crosses over $K_2$ \cite{Rolf}.

\begin{sect} {\bf Theorem.\/}
For any cocycle of $T_k$,
  where $k$
is any positive integer,
  and for any link $L$,
the state-sum $\Phi(L)$ is a function of
pairwise linking numbers.
\end{sect}
{\it Proof.\/}
Let  the elements of  $T_k$ be denoted by  $1, \ldots , k$, and let
$\phi= \prod_{i \neq j} \chi_{(i, j)} ^{w_{i,j}}
\in Z^2(T_k ; {\bf Z})$
(any cocycle can be written this way for some integers $w_{i,j}$).
For each coloring of $L$ by $T_k$, there is an ordered
partition
$A= \{ A_1, \ldots,  A_k \}$ of $\{ K_1 , \ldots K_n \}$
such that each component of $A_j$
is colored by $j \in T_k$ where $L= K_1 \cup \ldots \cup K_n$.
All ordered partitions of $\{ K_1 , \ldots , K_n \}$ are in one-to-one 
correspondence
to colorings by $T_k$.
Then the state-sum invariant $\Phi(L)$ with respect to this cocycle
is written as
$$ \sum_{A} \prod_{i , j } t^{ \lk (A_i, A_j) w_{i,j}  } $$
where
$A$ ranges over all
ordered
partitions of components.
  $\Box$

In particular, for $T_2$ and for links with small numbers of components,
we obtain the following formulas
by counting component-wise crossing numbers.

\begin{sect} {\bf Proposition.\/}
Take   $\phi=\chi_{(1,2)} \in Z^2(T_2 ; {\bf Z})$ to define a
cocycle invariant $ \Phi(L) = \Phi_{\phi}(L)$ for a link (or a knot) $L$.
\begin{enumerate}
\item
If $K$ is a knot, then $\Phi(K)=2$.
\item
If $L=K_1 \cup K_2$ is a $2$-component link, then
$$\Phi(L)= 2\left( 1+t^{\lk(L)} \right).$$
\item
If $L=K_1 \cup K_2 \cup K_3$ is a $3$-component link, then $$
\displaystyle \Phi(L)=  2\left(1 + \sum_{i,j=1,\;i<j} ^{3} t^{ \lk (L) -
\lk(K_i, K_j) } \right) .
\quad \Box  
$$
\end{enumerate}
\end{sect}

\vspace{5mm}

Next we study invariants with dihedral quandles.
We consider
  the dihedral quandle of four elements as $R_4 = \{a_1, a_2,
b_1, b_2 : a_i * a_j = a_i, b_i * b_j = b_i, a_i * b_j
= a_{i+1}, b_i * a_j = b_{i+1} \}$ where, in the subscripts,
$2+1$ is taken to be $1$.
Geometrically
$a_1$, $a_2$, $b_1$, and  $b_2$  represent the reflections
of a square about the horizontal axis, vertical axis,
the line $y=x$, the line $y=-x$, respectively.

\begin{figure}
\begin{center}
\mbox{
\epsfxsize=2in
\epsfbox{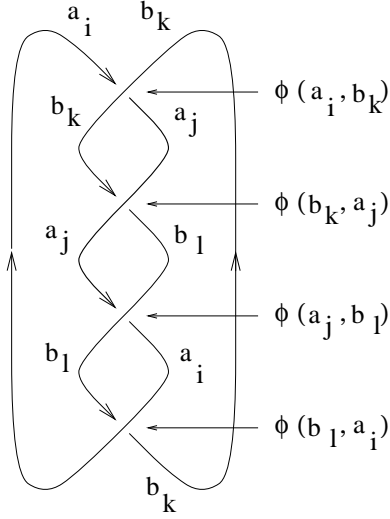}
}
\end{center}
\caption{Computations for $(4,2)$-torus link  }
\label{toruslink}
\end{figure}

\begin{sect} {\bf Example ($(4,2)$-torus link).\/} \label{torus42ex} {\rm
A computation of the state-sum invariant for the $(4,2)$-torus
link with $X=R_4$ is depicted in Fig.~\ref{toruslink}.
We assume $A= {\bf Z}= \langle t \rangle$.

First, we consider a
$2$-cocycle
$\phi = f_{(a_1, b_1)}= \chi_{(a_1, b_1)}  \chi_{(a_1, b_2)}$
(in multiplicative notation).
In Fig.~\ref{toruslink} a specific coloring
and the corresponding
weights are shown.
All possible colorings
are obtained as follows.
If only one quandle element is used, the coloring's
state-sum contribution
is trivial ($1$).  Since $R_4$ has 4 elements, there are 4 such
possibilities.  If one component is colored by $a_1$, and the other by $a_2$,
there are no crossings of weight $\phi(a_1,b_j)^{\pm 1}$ for $j \in\{1,2\}$,
  so these two colorings give trivial state-sum contributions as well.
Coloring one component by $b_1$ and the other by $b_2$ produces $2$ similar
cases.
When one component is colored by
the $a\/$s and  the other by the $b\/$s, the color contributes
  $t$ to the state-sum. There are $8$ such colorings,
one of which
is depicted in Fig.~\ref{toruslink}.
   Since these cases cover all possible colorings,
the state-sum is
$8+8t=8(1+t)$.

} \end{sect}

In general we have the following lemma on colorings by $R_4$,
which is proved by induction on the number of components $n$.

\begin{sect} {\bf Lemma.\/} \label{duallem}
Let $L=K_1\cup \ldots \cup K_n$ be a link such that any
pairwise linking number is even.
Then for any arcs $r_i$ of $K_i$, $i=1, \ldots, n$, and
for any color on $r_i$ ($i=1, \ldots, n$) by $R_4$,
there is a unique coloring of $L$ that extends the given coloring on $r_i$
  ($i=1, \ldots, n$).
In particular,
the number of colorings
is $4^n$.
\end{sect}

{}From the proof of Lemma~\ref{65},
any cohomology class  in $H^2(R_4; {\bf Z})$
is represented by a cocycle of the form
$\phi=\lambda_1^u  \lambda_2^v$
for some integers $u$ and $v$ where
$\lambda_1=\chi_{(a_1, b_1)} \chi_{(a_1, b_2)}$
and $\lambda_2=\chi_{(a_2, b_1)} \chi_{(a_2, b_2)}$.
Hence by Lemma~\ref{coblemma3}, all possible values of the state-sum invariant
with $R_4$ can be obtained by examining the cocycles of
the above form.

\begin{sect} {\bf Theorem.\/}
The state-sum invariant $\Phi(L)$  with respect to
the  cocycle $\phi=\lambda_1^u  \lambda_2^v$ of $R_4$, of
any $n$-component link $L=K_1 \cup \ldots \cup K_n$
such that any pairwise linking number is even,
is of the form
$$\Phi(L)= 2^n \left( \sum_{ A \in {\cal  P (K)}} t^{ (u+v) \lk(A, B) / 2}
\right) $$
where ${\cal K}=\{ K_1, \ldots, K_n \}$, $ {\cal  P (K)}$
denotes its power set, and
$B={\cal K} \setminus  A $.
\end{sect}

For example, for  $1$-,  $2$-, $3$-component links $L$, respectively,
\begin{eqnarray*}
\Phi(L) &=& 4 \\
\Phi(L) & = & 8\left(1+  t^{(u+v) \lk(L) /2 }  \right) \\
\Phi(L)  & =&  16\left(1 + \sum_{i,j=1,\;i<j} ^{3} t^{ (u+v) (\lk (L) -
  \lk(K_i, K_j)/2 )  } \right).
\end{eqnarray*}
{\it Proof.\/}
Write  $L= (\cup A)  \cup (\cup B) $ and take a color
that assigns $a_i$s to $A$ and $b_j$s to $B$.
Let $Y_i^{\pm}$ be the number of crossings
of a diagram of $L$ of sign $\pm$ where
the arc colored $a_i$ goes under an arc colored by
$b_1$ or $b_2$
and comes out with the color $a_{i+1}$.
Then with $\phi$, the state-sum contribution of
this color
is $t^{ u(Y_1^+ - Y_2^-) + v (Y_2^+ - Y_1^-) }$.
When we trace each  component $K_i$ of $A$,
  the colors alternate $a_1$ and $a_2$
at the crossings of the above types. Therefore
$$ Y_1^+ + Y_1^- =  Y_2^+ + Y_2^- , $$
which is equivalent to
$$  Y_1^+ -  Y_2^- =  Y_2^+ - Y_1^-, $$
so the contribution is written as
$$ t^{ u(Y_1^+ - Y_2^-) + v (Y_2^+ - Y_1^-) } = t^{ (u+v) ( Y_1^+ - Y_2^- 
)}. $$
On the other hand,  one computes
\begin{eqnarray*}
\lk(L) &=& ( Y_1^+  + Y_2^+ ) - ( Y_1^- + Y_2^-) \\
&=& 2( Y_1^+ - Y_2^-)
\end{eqnarray*}
and the result follows.
{}From the uniqueness in Lemma~\ref{duallem}, the number of such colorings
is (the number of colorings of components in $A$ by $a_i\/$s)
$\times$ (the number of colorings of components  in $B$ by $b_j\/$s)
  $=2^{|A|} \times 2^{|B|}=2^n$.
$\Box$

\vspace{5mm}

Recall that a
  map $f: X \rightarrow Y$ between two
quandles
$X, Y$ is called
a
(quandle) 
{\it homomorphism} if $f(a*b)=f(a)*f(b)$ for any $a, b \in X$.
A homomorphism
is called an {\it isomorphism } if
it is bijective.  
An isomorphism $f :X \rightarrow X$ is called an
{\it automorphism} (see \cite{Joyce}).

Next we consider invariants with general  dihedral quandles.
For a dihedral quandle $R_n$, use integers modulo $n$,
$R_n= \{ 0, 1, \ldots, n-1 \}$,  with the quandle operation
$i*j = 2j-i \pmod{n}$.
  Denote by $R_{2n}^0$
(respectively $R_{2n}^1$) the evens (resp. odds) of $R_{2n}.$
It is obvious that for any coloring of a link
of $L$ by $R_{2n}$,  each component
of $L$ is colored    either by
  $R_{2n}^0$ or by $R_{2n}^1$.

\begin{sect}{\bf  Lemma.\/}
\label{dihtriv}
If a knot $K$ has a nontrivial state-sum with $R_{2n}$
associated with a $2$-cocycle $\phi \in C^2(R_{2n}; {\bf Z})$,
then there is a $2$-cocycle $\phi' \in C^2(R_{n}; {\bf Z})$
such that $K$ has a nontrivial state-sum with $R_{n}$
associated with $\phi'$.
\end{sect}
{\it Proof.\/}
Let ${\cal C}$ be a coloring
such that
$\phi$ applied to the coloring
$\cal{C}$ of $K$ produces $t^q$-terms
for some $q \in {\bf Z}-0$; i.e.
$$  \prod_{\tau}  B( \tau, {\cal C} ) = t^q
\quad {\rm for~some ~} q \in {\bf Z}-0 . $$
Since $K$ is a knot, all colors
used in ${\cal C}$ are
either elements of $R^1_{2n}$
or elements of $R^0_{2n}$.
   For the first case, consider
an  isomorphism
$j_1 :R_{2n}^1  \rightarrow
R_n, $ $x \mapsto \frac{x-1}{2}$.
(The inverse is
$x \mapsto 2x+1$.)
The cocycle $\phi$ induces
a cocycle  $\phi' \in C^2(R_n; {\bf Z})$ by
$\phi'(x,y)= \phi(j_1^{-1}(x),j_1^{-1}(y))$.
Note that $\phi'$ has the same value on the color $j_1(\cal{C})$
as $\phi$
did on $\cal{C}$.  Thus, $K$ has a nontrivial
state-sum with $R_n$ associated with $\phi'$.
The other case is proved
similarly by use of $j_0: R_{2n}^0 \rightarrow R_n, $
$x \mapsto\frac{x}{2}$ (and the inverse
$x \mapsto 2x$). $\Box$

\begin{sect}{\bf Corollary.\/} All knots
have trivial state-sums with
any diherdal quandle $R_n$, $n \geq 2$,
associated with any $2$-cocycle $\phi \in C^2(R_n; {\bf Z})$. 
\end{sect}
{\it Proof.\/}  First note that $R_2$ is isomorphic to $T_2$, the
trivial two element quandle.  Then all knots in $R_{2}$ have
trivial state-sum.  Then,  by
induction using the above lemma, we see that any knot has the
trivial state-sum with  $R_{2^m}$ for all $m \geq 1$.
The dihedral quandles
with an odd number of elements have no 2-dimensional integral
quandle cohomology \cite{Greene}, \cite{cjks2}.
Hence the same argument, using the above lemma, gives the result. $\Box$

Next we consider  invariants with Alexander quandles.

\begin{sect} {\bf Lemma.\/} \label{LXquo}
  Let $n$ be a positive integer ($>1$)
and $h(T) \in {\bf Z}_n[T, T^{-1}]$.

{\rm (a)} If $n | h(1)$, then the map
$p: {\bf Z}_n[T, T^{-1}]/ (h(T))
\rightarrow  {\bf Z}_n$
defined by $f(T) \mapsto f(1)$ defines a surjective
homomorphism to the trivial quandle. Here $ {\bf Z}_n$ is given the quandle 
structure of $T_n$.

{\rm (b)} If $n|h(-1)$, then the  map $q: {\bf Z}_n[T, T^{-1}]/(h(T))
  \rightarrow  {\bf Z}_n$
defined by $f(T) \mapsto f(-1)$ defines a surjective
homomorphism to the dihedral quandle. Here $ {\bf Z}_n$ is given the 
quandle structure of $R_n$.
\end{sect}
{\it Proof.\/}
The operation
$a*b=Ta + (1-T)b$ on Alexander quandles become
$a*b=a$ for $T=1$ and $a*b=2b-a$ for $T=-1$.
$\Box$ 

\begin{sect} {\bf Remark.\/} 
\label{rem89} {\rm
(1) Observe that
part (b) corresponds to the existence of Fox colorings if and only if
$n$ divides the determinant of the knot.

(2) In the above lemma (a),
let $\phi_{i,j} \in Z^2(  {\bf Z}_n[T, T^{-1}]/(h(T)) ;  {\bf Z})$
be the pull-back cocycle
$p^{\sharp}\chi_{(i,j)}$
by the homomorphism $p$
of a cocycle $\chi_{(i,j)} \in  Z^2( T_n ; {\bf Z})$.
Then $\phi$ is written as $\phi=\prod \chi_{(f,g)}$ where the product
ranges over all $f, g \in {\bf Z}_n[T, T^{-1}]/(h(T))$
such that $f(1)=i$, $g(1)=j$ for $i, j \in  {\bf Z}_n$,
where $n|h(1)$.
We use this cocycle in the following theorem.
}\end{sect}

\begin{sect} {\bf Theorem.\/}
\label{noname}
For any positive integers $n, m >1$, there exists
a cocycle in $C^2({\bf Z}_n[T, T^{-1}]/(T^{2m}-1); {\bf Z})$ and 
a link $L$ whose cocycle invariant is non-trivial.

In particular, $H^2({\bf Z}_n[T, T^{-1}]/(T^{2m}-1); {\bf Z}) \neq 0$
for any $n,m>1$.
\end{sect}
{\it Proof.\/}
Let $L$ be the $(2mn, 2)$-torus link,
which is 
the closure of the braid $\sigma_1^{mn}$.
If the elements $a, b$ are assigned as colors to
the top of two strings of the braid, then after the $k\/$th
crossing, the colors assigned are
$[a,b]B^k$ where $B$ is the Burau matrix
$B= \left[ \begin{array}{cc} 0&T \\1&1-T \end{array} \right]$.
Each entry of the matrix $B^{2mn}-I$ is divisible by
the Alexander polynomial of $L$ (see \cite{K&P} for example).
The Alexander polynomial of $L$
is $\Delta=T^{2mn-1} - T^{2mn-2} + \ldots -1$
  (see for example \cite{Murasugi}).
With the relation $T^{2m}=1$ and with the coefficients in ${\bf Z}_n$,
$\Delta =0$.
Hence $B^{2mn}=I$ in ${\bf Z}_n[T, T^{-1}]/(T^{2m}-1)$,
and any pair $(a,b)$ gives a coloring
of $L$.
In particular, the pair $(0, 1)$ defines a coloring
and gives the term $T^{mn}$
with the cocycle
$\phi_{0,1}$
defined in Remark~\ref{rem89}. 
$\Box$

\begin{sect} {\bf Example.\/} {\rm
Let $L_{2n}$ be the $(2n, 2)$-torus link.
We use the Alexander quandle ${\bf Z}_3[T, T^{-1}]/(T^{2}-1)$.
Let $\phi=\phi_{0,1} \phi_{0,2}^2 \phi_{1,2}^3$
where $\phi_{i,j}$ are defined
in Remark~\ref{rem89}. 
By listing the colorings, one computes that
the cocycle invariant of $L_{2n}$ is
$27 + 18 (t^{3n} + t^{6n} + t^{9n})$
if $n=3m$, and
$9 + 6 (t^{3n} + t^{6n} + t^{9n})$
otherwise.
} \end{sect}

If we use the quandle
  $S_4= {\bf Z}_2[T, T^{-1}]/(T^2+T+1)$,
  then the trefoil
and  the figure 8 knot have non-trivial invariants using the
$2$-cocycle
$$\phi = t^{\chi_{(0,1)} +\chi_{(1,0)}+\chi_{(1+T,0)}+\chi_{(0,1+T)}
+ \chi_{(1,1+T)} +\chi_{(1+T,1)}} \in Z^2(S_4; {\bf Z}_2).$$

\begin{sect}{\bf Theorem.\/}
The state-sum invariant for the trefoil and the figure 8 knot with the
$2$-cocycle $\phi$ defined above is
$$4+12t.$$
\end{sect}
{\it Proof.} One can easily show that each knot can be colored in
  16 ways using this quandle.
The rest of the proof is a direct calculation.
$\Box$

\vspace{5mm}

Many
other knots also can be seen to have the polynomial
$4+12t$ and its
multiples as their invariants
using this cocycle.
It has also been computed that many knots in the knot table
have non-trivial values with a different quandle than $S_4$.

\section{
Triple Linking of Surfaces and Cocycle Invariants}
\label{sss2}
The linking number of a $2$-component classical
link $L=K_1 \cup K_2$ can be defined by counting the crossing number with signs
($\pm 1$) where the component $K_1$ crosses over $K_2$ (\cite{Rolf},
see also the preceding section).
This definition is generalized as follows to
linked surfaces.
Throughout
this section,
linked surfaces refer to
  oriented, multi-component, smoothly (or PL locally flatly) embedded
surfaces in $4$-space.

Recall from
Definition~\ref{triplepointsign} that
  the sign of a triple point is determined by comparing the ordered
triple of vectors normal to the top, middle, and bottom sheets to
  the right-handed
orientation  of $3$-space.
Let $F=K_1 \cup \cdots \cup K_n$ be a linked surface,
where
$K_i$, $i=1, \cdots, n$, are
components.

\begin{sect} {\bf Definition.\/} {\rm
Let $T_{\pm}(i,j,k)$ denote the number of
positive and negative, respectively, triple points
such that the top, middle, and bottom sheets are from
components
$K_i$, $K_j$, and $K_k$ respectively.
Such a triple point is called of type $(i,j,k)$.
Then define $T(i,j,k)=T_+(i,j,k) - T_-(i,j,k)$.
} \end{sect}

\begin{sect}{\bf Lemma.\/}
The numbers $T(i,j,k)$ are invariants of isotopy classes of $F$
if $i\neq j$ and $j \neq k$.
\end{sect}
{\it Proof.\/}
Consider the Roseman moves,
depicted in Fig.~\ref{rose},
that
are analogues of
the Reidemeister moves.
The invariance of $T(i,j,k)$ is proved by checking that they remain
unchanged
under these moves.

More specifically, there are three moves involving triple points:
(1) cancelation/creation of a pair of triple points (depicted in Fig~\ref{rose}
right top), (2) a  branch point passing through a sheet (left bottom),
and (3) the tetrahedral move, a move involving four planes
(right bottom).
In
  move (1), a pair of positive and negative triple points are
involved, so that the number
$T(i,j,k)$
remains
unchanged.
In
move (2), the triple point involved is of type $(i,i,j)$
or $(i,j,j)$ because the branch point connects two sheets in the
triple point, and these are the cases excluded in the Theorem.
The types of the various triple points remain
the same on either side of
move (3).
$\Box$

\vspace{5mm}

Thus these numbers are invariants of linked surfaces,
which we call {\it triple point linking invariants},
or simply {\it triple point invariants}.

Although we provided a diagrammatic definition and proof,
this invariant has been known in different contexts,
see \cite{Sanderson,Sanderson2,Koschorke} for example.

\begin{sect} {\bf Theorem.\/} \label{thm2}
For a linked oriented surface $L=K_1 \cup \cdots \cup K_n$ and for
any distinct $i, j \in \{1, \ldots n \}$, we have
$T(i,j,i)=0$.
\end{sect}
{\it Proof.\/}
Consider the double curves $D(i,j)$ where the over-sheet
is $K_i$ and the under-sheet is $K_j$.
Then $D(i,j)$ is a set of
immersed
closed curves.
The double curve is oriented in such a way that the ordered vectors of
the normals
$\vec{n}_i$ of $K_i$ and $\vec{n}_j$ of $K_j$
together with the direction $\vec{v}$
of $D(i,j)$ matches the orientation of ${\bf R}^3$.
Push each component of $D(i,j)$ off of $F$, to obtain a set of
closed oriented curves
$\gamma$ where the orientation is parallel to that of $D(i,j)$.
Then the intersection number, $\gamma \cap K_i$, counted with sign
is zero for a homological reason.
Such intersections occur near triple points of type
$(i,i,j)$ and $(i,j,i)$.
Near each
triple point of type $(i,i,j)$
a pair of   intersections occurs,
and they   occur in cancelling signed pairs.
However, near each triple point   of type $(i,j,i)$
a single intersection occurs,
and its sign matches the sign of
the triple point.
Therefore, $T(i,j,i)=0$.
$\Box$

\vspace{5mm}

\begin{sect} {\bf Theorem.\/}
For a linked oriented surface $L=K_1 \cup \cdots \cup K_n$ and for
any triple $(i,j,k)$ with $i,j,k
\in
\{ 1, \ldots , n\}$ where $i,j,$ and $k$ are distinct,
it holds that
$$ T(i,j,k)-T(i,k,j)+T(k,i,j)=0. $$
\end{sect}
{\it Proof.\/}
The same argument as in the proof of Theorem~\ref{thm2},
applied to
$D(i,j)$ and $\gamma \cap K_k$
gives
the  equality. Note that the middle term receives
a negative sign because of
the sign conventions of intersection and that of
triple points are opposite at these triple points. $\Box$

\vspace{5mm}

  The above conditions are equivalent to

\begin{sect} {\bf Corollary.\/} \label{solved}
For any three component  linked surface,
there exist integers $a$ and $b$ such that
$$ \begin{array}{llcrr}
T(1,2,3) &=& a& =& -T(3,2,1) \\
T(3,1,2) &=& b&   = & - T(2,1,3)  \\
T(2,3,1) & =& -(a+b) & = & -T(1,3,2)
\end{array} $$
\end{sect}

\begin{figure}
\begin{center}
\mbox{
\epsfxsize=2in
\epsfbox{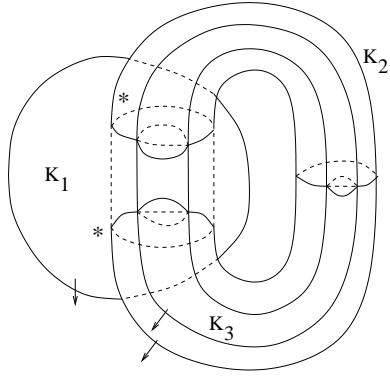}
}
\end{center}
\caption{A surface link with non-trivial triple linking  }
\label{nontrivial}
\end{figure}

\begin{sect} {\bf Theorem.\/} \label{tripthm}
For any integers $a$ and $b$, there exists a linked surface
$F=K_1 \cup K_2 \cup K_3$ such that $T(i,j,k)$
  $ (\{ i,j,k \} = \{ 1,2,3 \})$
satisfy
the conditions in Cor.~\ref{solved}.
\end{sect}
{\it Proof.\/}
Consider the linked surface
$F= K_1 \cup K_2 \cup K_3$ depicted in Fig.~\ref{nontrivial}.
Here $K_1$ is a sphere and $K_2$ and $K_3$ are tori.
Two tori $K_2$ and $K_3$ intersect along two parallel
double curves, such that $K_2$ is the over-sheet along one
of them, and $K_3$ is over-sheet along the other.
In other words, $K_2 \cup K_3$ is a spun Hopf link.
There are two intersections between $K_1$ and $K_2 \cup K_3$.
In the figure, at the bottom intersection $*$ is marked on the sphere,
indicating that $K_1$ is the top sheet over all sheets of $K_2 \cup K_3$,
and in the top intersection,
$*$
is placed on $K_2$ to
indicate that all sheets of  $K_2 \cup K_3$ are over $K_1$ (in other words
$K_1$ is the bottom).
Then the triple point invariants are computed as
$T(1,2,3)=1$, $T(1,3,2)=1$, $T(2,3,1)=-1$, and $T(3,2,1)=-1$.
This is the case where $a=1$ and $b=0$.
An example of a linked surface with  $a=0$ and $b=1$ is obtained
by switching the components, and the cases $a=-1$, $b=0$ and
$a=0$, $b=-1$ are obtained by changing orientations.
The general case is obtained by taking appropriate connected sum of
copies of these examples. $\Box$

\begin{sect} {\bf Theorem.\/}
Let $X=\{x,y,z \} $ be the trivial quandle of
three elements and $\theta \in Z^3(X; {\bf Z})$
be the cocycle
$\chi_{(x,y,z)}$ which is the characteristic function:
$$\chi_{(x,y,z)} 
(p,q,r)  
= \left\{ \begin{array}{lr} t & {\mbox{\rm if \ }}
  (p,q,r) 
= (x,y,z), \\ 1  & {\mbox{\rm otherwise.}}
\end{array} \right. $$
For a linked surface of three
connected
components
with triple point invariants
as given in Cor.~\ref{solved}, the
state-sum
invariant is
$$ t^a + t^{-a} + t^b + t^{-b} + t^{a+b} + t^{-a-b} + 21 .$$
\end{sect}
{\it Proof.\/}
There are 27 ways of coloring the link: For a given
connected component, each region of that component
has the same color as all the other regions of that component.

If a coloring uses fewer than three colors, then it contributes
the value 1 as a term in the state-sum. On the other hand,
for a given coloring ${\mathcal C}$, say
${\mathcal C}(K_i)=z$, ${\mathcal C}(K_j)=y$,
${\mathcal C}(K_k)=x$,
the Boltzmann weight of a triple point
is $t$ if and only if the triple point
is positive
and of type $(i,j,k)$.  The weight
is $t^{-1}$ if and only if the triple point is
negative of the same type. The weight is 1 otherwise.
So this coloring contributes  a term $t^{T(i,j,k)}$ to the
state-sum. $\Box$

\begin{sect} {\bf Remark.\/} {\rm
The same argument as above, together with Theorem~\ref{thm2},
gives
that when the two-element trivial quandle, $T_2$,
  is used,
the state-sum invariant 
assicated with 
any 3-cocycle
is trivial.
} \end{sect}

\section{ Surface Braids and Quandles}
\label{surfbraid}

In this section we give a method to obtain a
presentation of the quandle of a surface braid
described by a chart, which is used in order to
calculate the
state-sum
invariants of surfaces
in 4-space.

Let  $D^2$  and  $D$  be 2-disks and  $X_m$  a
fixed set of $m$  interior points of  $D^2$.
By  $pr_1 : D^2 \times D \to D^2$  and
$pr_2 : D^2 \times D \to D$, we mean the
projections to the first factor and to the
second factor, respectively.

\begin{sect}{\bf Definition.\/} {\rm
A {\em surface braid\/}
(\cite{Kam:ribbon,Rud})
  of degree  $m$ is a
compact, oriented surface  $S$  properly
embedded in  $D^2 \times D$  such that the
restriction of  $pr_2$  to  $S$  is a degree-$m$
simple branched covering map and  $\partial S
= X_m \times \partial D \subset D^2 \times \partial D$.
A degree-$m$ branched covering map $f : S \to D$
is {\em simple\/} if
$| f^{-1}(y) | = m$ or $m-1$ for $y \in D$.
In this case, the branch points are simple ($z \mapsto z^2$).

A surface braid  $S$  of degree  $m$  is extended
to a closed surface  $\widehat{S}$  in  $D^2 \times
S^2$  such that
$\widehat{S} \cap (D^2 \times D) = S$ and
$\widehat{S} \cap (D^2 \times \overline{D}) =
X_m \times \overline{D}$,
where  $S^2$  is the 2-sphere obtained from  $D^2$
by attaching a 2-disk  $\overline{D}$  along
the boundary.
     By identifying  $D^2 \times S^2$  with the
tubular neighborhood of a standard 2-sphere
in ${\bf R}^4$, we assume that  $\widehat{S}$  is a
closed oriented surface embedded in  ${\bf R}^4$.
We call it the {\em closure\/}  of $S$  in ${\bf R}^4$.
It is proved in
\cite{Kam:top}
that every closed oriented surface embedded in ${\bf R}^4$ is
ambient isotopic to the closure of a surface braid.

Two surface braids  $S$  and  $S'$  in $D^2 \times D$
are said to be {\em equivalent\/}  if there is an
isotopy  $\{ h_t \}$  of  $D^2 \times D$  such that
\begin{enumerate}
\item
$h_0 = {\rm id}$, $h_1(S) = S'$,
\item
for each $t \in [0,1]$, $h_t$ is fiber-preserving;
that is, there is a homeomorphism
$\underline{h}_t : D \to D$ with
$\underline{h}_t \circ pr_2 = pr_2 \circ h_t$, and
\item
for each $t \in [0,1]$, $h_t |_{D^2 \times \partial D}
= {\rm id}$.
\end{enumerate}

Let  $C_m$  be the configuration space of unordered
$m$  interior points of  $D^2$.  We identify the
fundamental group  $\pi_1(C_m, X_m)$  of
$C_m$  with base point  $X_m$  with the
braid group  $B_m$
on $m$ strings.
     Let  $S$
denote
a surface braid and $\Sigma(S) \subset D$
the
branch point set of the branched covering map
$S \to D$.
For a path
$ a : [0,1] \to D  \setminus
\Sigma(S)$,  we define a path
$$ \rho_S (a) : [0,1] \to C_m  $$
by
$$ \rho_S (a) (t) =
pr_1( S \cap (D^2 \times \{ a(t) \})). $$
     If  $pr_1( S \cap (D^2 \times \{ a(0) \})) =
pr_1( S \cap (D^2 \times \{ a(1) \})) = X_m$,
then the path  $\rho_S (a)$ represents an element
of  $\pi_1(C_m, X_m) = B_m$.
     Take a point  $y_0$ in $\partial D$.
The {\em braid monodromy\/}  of  $S$  is the
homomorphism
$$ \rho_S : \pi_1( D  \setminus \Sigma(S), y_0 ) \to B_m $$
such that  $\rho_S ([a]) = [ \rho_S (a) ]$
for any loop  $a$  in  $D \setminus \Sigma(S)$
with
base point  $y_0$.
}\end{sect}

Let $\Sigma (S) = \{y_1, \ldots, y_n \}$. Take a regular neighborhood
$N(\Sigma (S)) = N(y_1)\cup \cdots \cup N(y_n)$ in $D$.
A {\em Hurwitz arc system\/}
${\cal A} = ( \alpha_1, \dots, \alpha_n )$  for $\Sigma(S)$
is an $n$-tuple of simple arcs in
$E(\Sigma (S)) = {\mbox{\rm Cl}}(D \setminus N(\Sigma (S)) )$
(where Cl denotes the closure)
such that each
$\alpha_i$ starts from a point of $\partial N (y_i)$ and ends at $y_0,$ and
$\alpha_i \cap \alpha_j = \{ y_0 \}$  for
$ i \neq j $,  and $\alpha_1, \ldots, \alpha_n$
appear in this order around  $y_0$.

Let $\eta_i$ ($i= 1, \dots, n$)  be
the
loop
$\alpha_i^{-1} \cdot \partial N(y_i)\cdot \alpha_i$
in $D \setminus  \Sigma(S)$
with base point  $y_0$
which goes along $\alpha_i$, turns
along $\partial N(y_i)$
in the positive direction,
and returns along  $\alpha_i$.

\begin{sect}{\bf Definition.\/} {\rm
The {\em braid system\/}  of $S$  associated with
${\cal A}$ is an  $n$-tuple of  $m$-braids
$$
( \rho_S ([\eta_1]), \rho_S ([\eta_2]), \dots,
\rho_S ([\eta_n]) ). $$

Each element of a braid system is a conjugate
of a standard generator  $\sigma_i$  of $B_m$
or its inverse.  The braid system of a surface
braid of degree $m$ is written as
$$
( w_1^{-1} s_1^{\epsilon_1} w_1,
w_2^{-1} s_2^{\epsilon_2} w_2,
\dots,
w_n^{-1} s_n^{\epsilon_n} w_n ), $$
where $n$ is the number of branch points,
$w_1, \dots, w_n$ are $m$-braids,
$s_1, \dots, s_n
\in \{ \sigma_1, \dots, \sigma_{m-1} \}$ and
$\epsilon_1, \dots, \epsilon_n \in
\{ +1, -1 \}$.
}\end{sect}

\begin{sect}{\bf Definition.\/} {\rm
An {\it $m$-chart} \cite{Kam:ribbon}
is
an
oriented, labelled graph $\Gamma$ in $D$, which may
be empty or have closed edges without vertices
(which are called {\em hoops\/}),
  satisfying the following conditions:
\begin{enumerate}
\item
Every vertex has degree one, four or six.
\item
The labels of edges are in $\{ 1, 2, \dots, m-1 \}$.
\item
For each degree-six vertex, three consective edges are
oriented inward and the other three are outward,
and these six edges are labelled  $i$  and  $i+1$
alternately for some  $i$.
\item
For each degree-four vertex, diagonal edges have the
same label and are oriented coherently, and
the labels  $i$  and  $j$  of the diagonals satisfy
$| i-j | > 1$.
\end{enumerate}
We call
a
degree 1 (resp. degree 6) vertex a
{\em black\/} (resp. {\em white\/}) vertex.
A degree 4 vertex is called a {\it crossing point of
the chart}.

We say that a path  $\alpha: [0,1] \to D$ is
{\em in general position with respect to $\Gamma$\/}
if it avoids the vertices of $\Gamma$ and every
intersection of  $\alpha$  and  $\Gamma$  is a
transverse double point.   If  $p$  is an intersection
of  $\alpha$  and  an  edge of $\Gamma$  labelled  $i$
and if the edge is oriented from right to left
(resp. from left to right), then assign the
intersection  $p$ a letter  $\sigma_i$
(resp. $\sigma_i^{-1}$).
Read the letters assigned the intersections of
$\alpha$ and $\Gamma$  along $\alpha$  and we have a word
$$
\sigma_{i_1}^{\epsilon_1} \sigma_{i_2}^{\epsilon_2}
\dots \sigma_{i_s}^{\epsilon_s}
$$
in the braid generator.  We call this the
{\em intersection braid word\/}  of $\alpha$
with respect to $\Gamma$, and denote it by
$w_\Gamma (\alpha)$.

For an $m$-chart $\Gamma$,
a {\it surface braid described by $\Gamma$} means
a surface braid
$S$ of degree $m$ satisfying the following conditions:
\begin{enumerate}
\item
For a regular neighborhood  $N(\Gamma)$  of $\Gamma$
in $D$ and
for any  $y \in {\mbox{\rm Cl}} ( D \setminus N(\Gamma))$,
the projection $pr_1$ satisfies the condition:
  $pr_1( S \cap (D^2 \times \{ y \}) ) = X_m$,
where $X_m$ denotes the $m$ fixed interior points of $D^2$.

\item
The branch point set of  $S$
corresponds to the set
of the black vertices of  $\Gamma$.
\item
For a path  $\alpha: [0,1] \to D$  which is in general
position with respect to $\Gamma$  and  $\alpha(0)$,
$\alpha(1)$ are in
${\mbox{\rm Cl}} ( D \setminus N(\Gamma))$, the $m$-braid determined by 
$\rho_S (\alpha)$  is
the  $m$-braid presented by the intersection braid
word  $w_\Gamma (a)$.
\end{enumerate}
}\end{sect}

\begin{sect} {\bf Proposition \cite{Kam:ribbon}.\/}
\begin{enumerate}
\item
For any $m$-chart $\Gamma$, there is a
unique (up to equivalence) 
surface braid
described by $\Gamma$.
\item
For any surface braid  $S$ of degree  $m$, there is an
$m$-chart  $\Gamma$  such that  $S$  is equivalent to
a surface braid described by  $\Gamma$.
\end{enumerate}
\end{sect}

\begin{figure}
\begin{center}
\mbox{
\epsfxsize=3in
\epsfbox{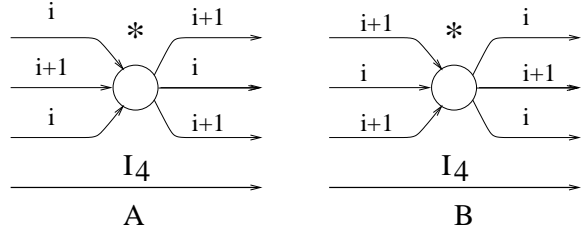}
} \end{center}
\caption{The distinguished region of
  a white vertex
}
\label{region}
\end{figure}

\begin{figure}
\begin{center}
\mbox{
\epsfxsize=2in
\epsfbox{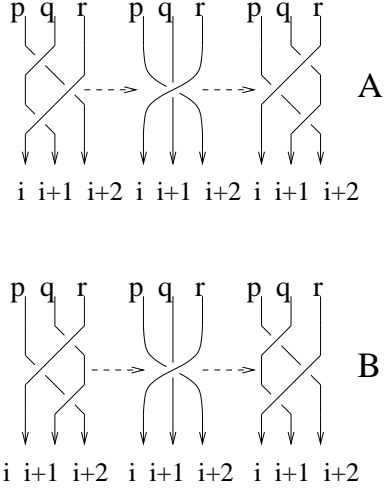}
} \end{center}
\caption{Quandle labels  near a white vertex
}
\label{quandlelabels}
\end{figure}

Let  $S$  be a surface braid described by
$\Gamma$.  Identify  $D^2$  with  $I_1\times I_2$
and  $D$  with  $I_3\times I_4$, where  $I_i$ ($i=1, \dots, 4$)
are
intervals.
For each  $t \in I_4$, put  $b_t = S \cap (D^2 \times I_3 \times \{t\})$.
Then  $\{ b_t | t \in I_4 \}$  is a continuous sequence of
$m$-braids with a finite number of exceptions that are singular
$m$-braids.   Modifying  $\Gamma$  by an ambient isotopy of  $D$,
we may assume that every white vertex  $W$  looks like one of
the Fig.~\ref{region}
with respect to the bi-parametrization
$D \cong  I_3 \times I_4$.  Then the sequence  $\{ b_t \}$ looks like
the
motion pictures
in Fig.~\ref{quandlelabels}
around the white vertex.

\begin{figure}
\begin{center}
\mbox{
\epsfxsize=3in
\epsfbox{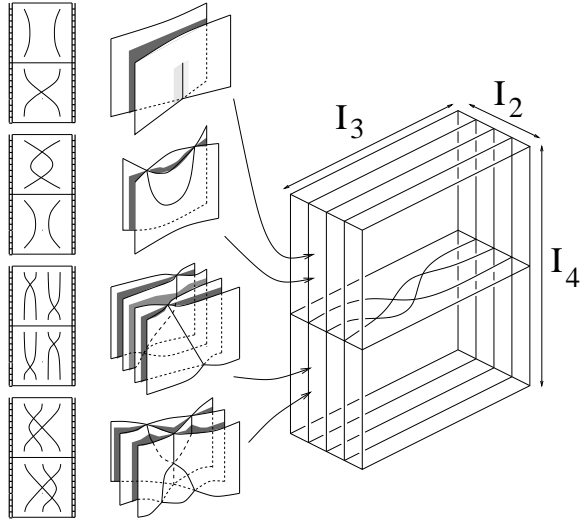}
} \end{center}
\caption{Projections and charts
}
\label{ebcproj}
\end{figure}

We assume that each $b_t$  is illustrated as a diagram
with respect to the projection  $I_1 \times I_2 \times I_3
\to I_2 \times I_3$.   Then under the projection
$I_1 \times I_2 \times I_3 \times I_4
\to I_2 \times I_3 \times I_4$, the image of $S$ has a
triple
point corresponding 
to a 
white vertex.
We define the {\em sign \/} $\epsilon (W)$ of a white vertex, $W$,
by  $+1$  (resp.  $-1$) if it is as (A) (resp. as (B)) in
Fig.~\ref{region}
  so that the corresponding triple point has
sign  $+1$   (resp.  $-1$) in the broken surface diagram of $S$.
(In general, the singularity set of the image of $S$ by the
projection  $I_1 \times I_2 \times I_3 \times I_4
\to I_2 \times I_3 \times I_4$  is identified
naturally
with the chart $\Gamma$
in the sense of
\cite{CS:book, Kam:Nato} 
The white vertices are in one-to-one correspondence to
the triple points and the black vertices
are to
the branch points.  Figure~\ref{ebcproj} shows the
relationship schematically, see
\cite{CS:book, Kam:Nato} 
  for details.)

When a surface braid is described by a chart, the
braid system is easily obtained as follows:

\begin{sect} {\bf Lemma.\/}
Let $S$  be a surface braid described by a chart  $\Gamma$,
and ${\cal A} = (\alpha_1, \dots, \alpha_n)$ a Hurwitz
arc system for $\Sigma(S)$ such that each $\alpha_i$ is
in general position with respect to $\Gamma$.
The braid system of  $S$
associated with  ${\cal A}$  is given by
$$
(w_\Gamma (\eta_1), w_\Gamma (\eta_2), \dots,
w_\Gamma (\eta_n)), $$
where $\eta_1, \dots, \eta_n$ are loops in
$D \setminus
\Sigma(S)$ associated with ${\cal A}$ as before.
\end{sect}
{\it Proof.\/}
By definition of  $S$,
$\rho_S ([\eta_i]) = w_\Gamma (\eta_i)$ for $i=1, \dots, n$.
$\Box$

\begin{figure}
\begin{center}
\mbox{
\epsfxsize=3in
\epsfbox{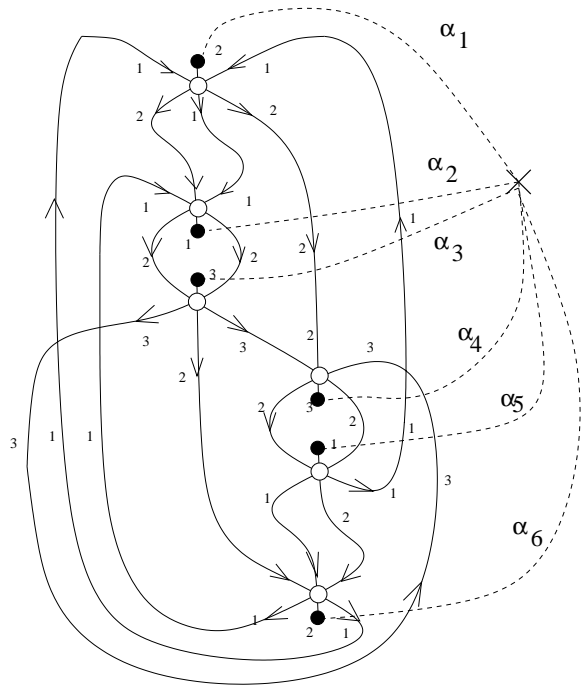}
} \end{center}
\caption{A Hurwitz system of a chart
}
\label{twisttre1holo}
\end{figure}

\begin{sect} {\bf Example.\/}
\label{anexample}
{\rm
Let  $\Gamma$ be
the $4$-chart as in
Fig.~\ref{twisttre1holo}
and  $S$  a surface braid of degree $4$ described by $\Gamma$.
It is known that this chart represents the $2$-twist spun trefoil knot
\cite{Kam:ribbon}.
Recall that every black vertex stands for a branch point
of $S \to D$.
Let ${\cal A} = (\alpha_1, \dots, \alpha_n)$ be a Hurwitz
arc system for $\Sigma(S)$ illustrated in the figure,
where $\alpha_1, \dots, \alpha_6$ are drawn as dotted arcs.
The braid system $(
w_1^{-1} \sigma_{k_1}^{\epsilon_1} w_1,
w_2^{-1} \sigma_{k_2}^{\epsilon_2} w_2,
\dots,
w_6^{-1} \sigma_{k_6}^{\epsilon_6} w_6)$
  of  $S$ is given by
$$
\begin{array}{ll}

w_1 = 1, \quad &
\sigma_{k_1}^{\epsilon_1} = \sigma_2^{-1}, \\

w_2 = \sigma_2^{-2}\sigma_1, \quad &
\sigma_{k_2}^{\epsilon_2} = \sigma_1, \\

w_3 = \sigma_2^{-2}\sigma_1, \quad &
\sigma_{k_3}^{\epsilon_3} = \sigma_3^{-1}, \\

w_4 = \sigma_2^{-1}\sigma_1\sigma_3, \quad &
\sigma_{k_4}^{\epsilon_4} = \sigma_3, \\

w_5 = \sigma_2^{-1}\sigma_1\sigma_3, \quad &
\sigma_{k_5}^{\epsilon_5} = \sigma_1^{-1}, \\

w_6 = \sigma_1^{-1}\sigma_3, \quad &
\sigma_{k_6}^{\epsilon_6} = \sigma_2.

\end{array}
$$
} \end{sect}

\begin{sect}{\bf Definition \cite{FR,Joyce}. \/}
{\rm
Let $M$
be an oriented $(n+2)$-manifold,
and $L$ be
an oriented $n$-submanifold of $M$ with a tubular neighborhood $N(L)$ in $M$.
Take a point $z \in E(L)= {\mbox{\rm Cl}}(M\setminus N(L)).$
Consider 
the set of paths $\alpha: [0,1] \rightarrow E(L)$ such that
there is a meridian disk, say $\Delta_\alpha$, of $L$ with
$\alpha(0) \in \partial \Delta_\alpha$ and $\alpha(1) = z.$ Let
$Q(M,L,z)$
be the set of homotopy classes of paths $\alpha$.
Define a binary operation $*$ on $Q(M,L,z)$
by
$$ [\alpha ] * [\beta]= [ \alpha \cdot \beta^{-1} \cdot  \partial 
\Delta_\beta \cdot \beta]$$
where $\Delta_\beta$ is an (oriented) meridian disk with $\beta(0) \in 
\partial \Delta_\beta$.
Then  $Q(  M,L,z )$  with  $\ast$  is a quandle,
which is called the {\em quandle of $(M,L )$\/},
or the {\em quandle of $L$\/}, with base point $z$,
and denote by $Q( M,L,z)$
(or $Q( M,L )$, $Q( L  )$, {\it etc.}).
}\end{sect}

\begin{figure}
\begin{center}
\mbox{
\epsfxsize=3in
\epsfbox{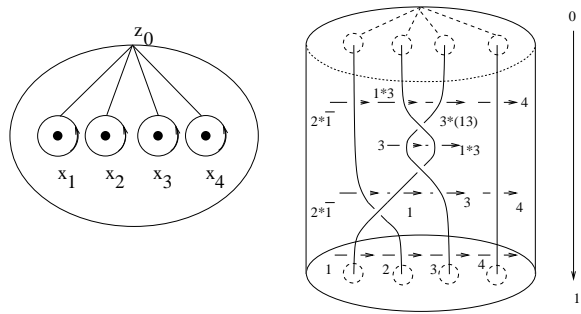}
} \end{center}
\caption{Generators of the free quandle
}
\label{braidquandle}
\end{figure}

\begin{sect}{\bf Example.\ }
{\rm
Let  $b$ be an $m$-braid, and  let
$f_b :
( D^2, X_m ) \to ( D^2, X_m ) $ be
an 
Artin homeomorphism associated with $b$.  We
denote by  $Q(b)$ the quandle isomorphism
$$(f_b)_*:
Q( D^2, X_m, z_0 ) \to Q( D^2, X_m, z_0 ) $$
induced from the Artin homeomorphism $f_b$.
    We usually identify
$Q( D^2, X_m, z_0 )$ with the free quandle
$F_Q \langle x_1, \dots, x_m \rangle$ generated by
$x_1, \dots, x_m$ as in
Fig.~\ref{braidquandle}
and regard  $Q(b)$  as a quandle automorphism of
the free quandle
$F_Q \langle x_1, \dots, x_m \rangle$.

The
quandle automorphism  $Q(b)$
is
interpreted as follows:
Let $\alpha$  be
the  path in $D^2 \times [0,1]$ defined by
$\alpha (t) = ( z_0, t )$.  We have
a quandle
isomorphism
$$
\alpha_\ast :
Q( D^2 \times [0,1], b,
z_0 \times \{ 0 \})  \to
Q( D^2 \times [0,1],  b,
z_0 \times \{ 1 \}) $$
such that
$\alpha_\ast ([\beta]) = [ \beta \cdot \alpha]$ .
Since the inclusion-induced
quandle
homomorphisms
$$
(i_j)_\ast :
Q( D^2 \times \{ j \}, X_m \times \{j \},
z_0 \times \{ j \})  \to
Q( D^2 \times [0,1], b,
z_0 \times \{ j \})
\quad {(  j= 0, 1  )}$$
are isomorphisms, we have an isomorphism
$$
(i_1)_\ast^{-1} \circ \alpha_\ast \circ (i_0)_\ast :
Q( D^2  \times \{ 0 \}, X_m \times \{ 0 \},
z_0 \times \{ 0 \})  \to
Q( D^2  \times \{ 1 \},X_m \times \{ 1 \},
z_0 \times \{ 1 \}). $$
Identifying
$Q( (D^2   \times \{ j \},X_m \times \{ j \},
z_0 \times \{ j  \}) $  ($j=0,1$) with
$Q( D^2, X_m, z_0 )$ via the
projection $D^2 \times [0,1] \to D^2$,
we have an automorphism of
  $Q( D^2, X_m, z_0 )$.
This is
$Q(b)$.

For example, if $b = \sigma_2^{-2}\sigma_1
\in B_4$, then the quandle
$Q( D^2, X_4, z_0 )$ is freely generated by
$x_1, \dots, x_4$ illustrated as in
Fig.~\ref{braidquandle}
and the quandle isomorphism  $Q(b)$ maps the
generators as follows:
\begin{eqnarray*}
Q(b) (x_1)  &=& x_2 \ast x_1^{-1}, \\
Q(b) (x_2)  &=& x_1 \ast x_3, \\
Q(b) (x_3)  &=& x_3 \ast (x_1 x_3), \\
Q(b) (x_4)  &=& x_4. \\
\end{eqnarray*}

In the above table and in the sequel,
we are mimicing the notation in \cite{FR}.
So, $a*(bc)$ is defined to be $(a*b)*c$;
  the element $a*(b^{-1})$ is the unique element
$c$ such that $a=c*b$; and generally a product $a*w$ where $w$
  is  a word on the free group generated by the quandle can be
  interpreted inductively. For example,
$x_1*(x_3 x_4 x_3^{-1}) = ((x_1 * x_3)*x_4)*x_3^{-1}$.
See also Fig.~\ref{braidquandle}.
}\end{sect}

\begin{sect} {\bf Lemma.\/}
\label{present}
Let $( b_1, \dots, b_n )$ be a braid system of
the surface braid
$S$,
then the quandle $Q(S)$ has a presentation
whose generators are
  $x_1, \dots, x_m$ and the relations are
$$
Q(w_i) (x_{k_i}) = Q(w_i) (x_{k_i + 1})  \quad
(i= 1, \dots, n), $$
where  $b_i = w_i^{-1} \sigma_{k_i}^{\epsilon_i} w_i$.
\end{sect}
{\it Proof.\/}
In \cite{Rud} and
\cite{Kam:top} it is shown that
the
fundamental
group
$\pi_1( D^2 \times D \setminus S, z_0 \times y_0 )$ is
generated by $m$ positive meridional elements
$x_1, \dots, x_m$ with defining relations
$$
(f_{w_i})_\ast (x_{k_i}) = (f_{w_i})_\ast (x_{k_i + 1})  \quad
(i= 1, \dots, n), $$
where
$(f_{w_i})_\ast : \pi_1( D^2
\setminus
X_m, z_0 ) \to
\pi_1( D^2
\setminus
X_m, z_0 )$ is the automorphism induced
from the Artin homeomorphism $f_{w_i} : ( D^2, X_m, z_0 ) \to
\pi_1( D^2, X_m, z_0 )$ associated with the braid
$w_i$.
In \cite{FR,Joyce}, presentations of quandles of
codimension $2$ embeddings in Euclidean spaces were given, that are similar to
Wirtinger presentations of fundamental groups.
Thus a similar argument as above  gives the presentation of $Q(S).$
  $\Box$

\begin{sect} {\bf Example.\/} {\rm
Let  $S$  be
the
surface braid of degree  $4$  described by
a $4$-chart $\Gamma$
  in Fig.~\ref{twisttre1holo}.
For a Hurwitz arc system
${\cal A} = ( \alpha_1, \dots, \alpha_6 )$ as in Fig.~\ref{twisttre1holo},
the braid system
$( w_1^{-1} \sigma_{k_1}^{\epsilon_1} w_1, \dots,
w_6^{-1} \sigma_{k_6}^{\epsilon_6} w_6 )$
of  $S$  is given as
in Example~\ref{anexample}.
%%%%%%%%%%%%%%don't delete for safety:%%%%%
% follows:
%$$
%\begin{array}{ll}%
%
%w_1 = 1, \quad &
%\sigma_{k_1}^{\epsilon_1} = \sigma_2^{-1}, \\
%
%w_2 = \sigma_2^{-2}\sigma_1, \quad &
%\sigma_{k_2}^{\epsilon_2} = \sigma_1, \\
%
%w_3 = \sigma_2^{-2}\sigma_1, \quad &
%\sigma_{k_3}^{\epsilon_3} = \sigma_3^{-1}, \\
%
%w_4 = \sigma_2^{-1}\sigma_1\sigma_3, \quad &
%\sigma_{k_4}^{\epsilon_4} = \sigma_3, \\
%
%w_5 = \sigma_2^{-1}\sigma_1\sigma_3, \quad &
%\sigma_{k_5}^{\epsilon_5} = \sigma_1^{-1}, \\
%
%w_6 = \sigma_1^{-1}\sigma_3, \quad &
%\sigma_{k_6}^{\epsilon_6} = \sigma_2.
%
%\end{array}
%$$
%%%%%%%%%%%%%%%%%%%%%%%%%%%%%%%
The quandle automorphisms
$Q(1)$,
$Q(\sigma_2^{-2}\sigma_1)$,
$Q(\sigma_2^{-1}\sigma_1\sigma_3)$, and
$Q(\sigma_1^{-1}\sigma_3)$
of
$F_Q \langle x_1, \dots, x_m \rangle$  map the generators as follows.

\begin{eqnarray*}
Q(1) &: &
x_1 \mapsto x_1, \quad
x_2 \mapsto x_2, \quad
x_3 \mapsto x_3, \quad
x_4 \mapsto x_4, \\
Q(\sigma_2^{-2}\sigma_1) &: &
x_1 \mapsto x_2 \ast x_1^{-1}, \quad
x_2 \mapsto x_1 \ast x_3, \quad
x_3 \mapsto x_3 \ast (x_1 x_3), \quad
x_4 \mapsto x_4, \\
Q(\sigma_2^{-1}\sigma_1\sigma_3) &: &
x_1 \mapsto x_2 \ast x_1^{-1}, \quad
x_2 \mapsto x_4 \ast x_3^{-1}, \quad
x_3 \mapsto x_1 \ast (x_3 x_4 x_3^{-1}), \quad
x_4 \mapsto x_3, \\
Q(\sigma_1^{-1}\sigma_3) &: &
x_1 \mapsto x_2, \quad
x_2 \mapsto x_1 \ast x_2, \quad
x_3 \mapsto x_4 \ast x_3^{-1}, \quad
x_4 \mapsto x_3.
\end{eqnarray*}

Hence the defining relations
$Q(w_i)( x_{k_i} ) = Q(w_i)( x_{k_i +1} )$ ($i=1,\dots, 6$) of
$Q(S)$  are
\begin{eqnarray*}
x_2  &=&
x_3, \\
x_2 \ast x_1^{-1}  &=&
x_1 \ast x_3,  \\
x_3 \ast (x_1 x_3)  &=&
x_4, \\
x_1 \ast (x_3 x_4 x_3^{-1})  &=&
x_3, \\
x_2 \ast x_1^{-1}  &=&
x_4 \ast x_3^{-1}, \\
x_1 \ast x_2  &=&
x_4 \ast x_3^{-1}.
\end{eqnarray*}

Thus the quandle $Q(S)$ is
\begin{eqnarray*}
\langle x_1, \dots, x_4 |
& & x_2 = x_1 \ast (x_2 x_1), \\
& & x_2 = x_2 \ast (x_1^2), \\
& & x_3 = x_2, \\
& & x_4 = x_1
\rangle \\
=
\langle x_1, x_2 |
& & x_2 = x_1 \ast (x_2 x_1), \\
& & x_2 = x_2 \ast (x_1^2)
\rangle.
\end{eqnarray*}

} \end{sect}

\section{Cocycle Invariants and Braid Charts}
\label{themain}

In this section we introduce a method to calculate the
state-sum
invariant of a surface braid described by a chart.
The
state-sum
invariant
of a surface braid
  coincides with the
state-sum
invariant of
its closure in  ${\bf R}^4$.

Let  $S$  be a surface braid of degree  $m$  described by an
$m$-chart  $\Gamma$.
The region of $D
\setminus
\Gamma$
assigned
the asterisk in Fig.~\ref{region}
is called the {\em distinguished region\/} for
a
white vertex  $W$.
Let  $y$  be a point of this region.  Since  $S$ is a surface
braid described by $\Gamma$, we may assume that
$pr_1( S \cap (D^2 \times \{ y \}) ) = X_m$.   Then
$Q_y = Q( D^2 \times \{ y \}, X_m \times \{ y \}, z_0 \times \{ y \})$
is identified with
$Q( D^2, X_m, z_0) = F_Q \langle x_1, \dots, x_m \rangle$
via the projection $pr_1$.   Take a path  $\beta : [0,1] \to
D \setminus \Sigma(S)$  with  $\beta (0) = y$  and   $\beta (1) = y_0$.
The $m$-braid $\rho_S(\beta)$
induces an isomorphism
$$ Q(\rho_S (\beta) ) :
F_Q \langle x_1, \dots, x_m \rangle = Q_y \to
Q_{y_0} = F_Q \langle x_1, \dots, x_m \rangle. $$
Recall that the generators  $x_1, \dots, x_m$  of the
quandle  $Q(S) = Q( D^2 \times D, S, z_0 \times y_0 )$  in
Lemma~\ref{present}
are the images  $i_\ast(x_1), \dots, i_\ast(x_m)$  of
$x_1, \dots, x_m \in Q_{y_0}$  by the inclusion-induced
homomorphism  $i_\ast : Q_{y_0} \to Q(S)$.

Let $p, q, r$  be the elements of  $Q(S)$  corresponding
to the
three sheets in the broken surface diagram of $S$ as in
the motion pictures depicted in
Fig.~\ref{quandlelabels}.
Then
\begin{eqnarray*}
p &= & i_\ast \circ Q(\rho_S (\beta) ) (x_i),  \\
q &= & i_\ast \circ Q(\rho_S (\beta) ) (x_{i+1}), \quad {\rm and}\\
r &= & i_\ast \circ Q(\rho_S (\beta) ) (x_{i+2}).
\end{eqnarray*}
We call  $(p,q,r)$  the {\em quandle triple\/} for the
white vertex $W$ or for the triple point corresponding to $W$.

Let  $\theta \in Z^3(X; A)$
be a 3-cocycle of a finite
quandle  $X$  with coefficient group  $A$  and let
$c : Q(S) =  Q( D^2 \times D, S, z_0 \times y_0 ) \to X$
be a
homomorphism (a coloring).
We define the {\em Boltzmann weight\/} on
a white vertex $W$ by
$$ \theta ( c(p), c(q), c(r) )^{\epsilon(W)}, $$
where  $\epsilon (W)$  is the sign of  $W$, and  $(p,q,r)$  is
the quandle triple for  $W$,  and put
$$ \Phi_\theta (\Gamma) =
\sum_c \prod_W \theta ( c(p), c(q), c(r) )^{\epsilon (W)}, $$
where $W$ runs over all white vertex of the chart $\Gamma$
and  $c$  runs over all possible coloring from
$Q(S)$ to $X$.

\begin{sect}{\bf Lemma.\/}
Let  $S$  be a surface braid described by a chart $\Gamma$  and
let $\widehat{S}$  be the closure of  $S$  in
${\bf R}^4$.
Then
$$ \Phi_\theta (\Gamma) = \Phi_\theta (\widehat{S}). $$
\end{sect}
{\it Proof.\/}
Recall the situation depicted in Fig.~\ref{ebcproj}.
Consider a broken surface diagram of  $S$  by the projection
$I_1 \times I_2 \times I_3 \times I_4 \to
I_2 \times I_3 \times I_4$.  The broken surface diagram of
$\widehat{S}$  is obtained from the diagram of  $S$ by
attaching  $m$  disks outside of  $I_2 \times I_3 \times I_4$
trivially.  So there is a one-to-one correspondence between
the colorings of them.  Every white vertex
corresponds
to a triple point, and
the Boltzmann weight of a
white vertex
is
the same with
that
of the triple point.
(In fact, we defined
it
to be so.)  Hence
$ \Phi_\theta (\Gamma) = \Phi_\theta (\widehat{S}). $
$\Box$

\vspace{5mm}

\begin{figure}
\begin{center}
\mbox{
\epsfxsize=3in
\epsfbox{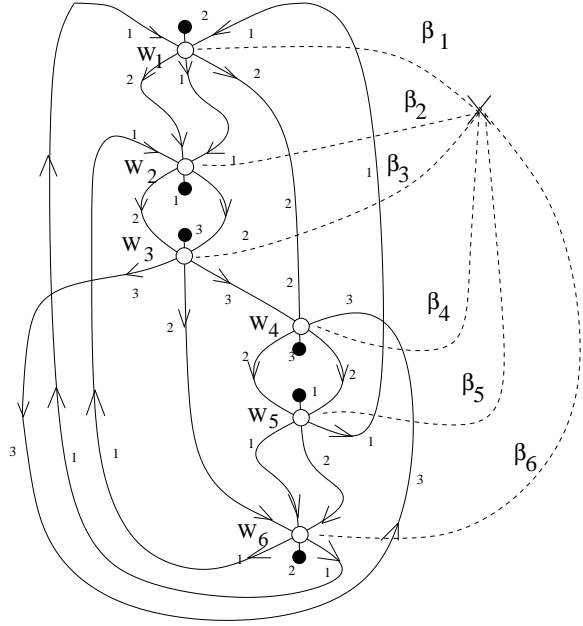}
} \end{center}
\caption{Computing cocycle invariants using a chart
}
\label{twisttre1cocy}
\end{figure}

\begin{sect} {\bf Theorem.\/}
\label{generalform}
Let  $F$  be the $2$-twist spun
trefoil and  $\theta$ a
quandle $3$-cocycle of a finite quandle  $X$ with
coefficient group  $A$.  The
state-sum
invariant
$\Phi_\theta (F)$ is
\begin{eqnarray*}
\sum_{y_1,y_2}
&&
\theta (y_1 \ast y_2, y_1, y_2)
\theta (y_1 \ast y_2, y_2, y_1 \ast y_2)
\theta (y_2, y_1 \ast y_2, y_1) \\
&&
\theta (y_1, y_1 \ast y_2, y_2)^{-1}
\theta (y_1 \ast y_2, y_1, y_1 \ast y_2)^{-1}
\theta (y_1, y_2, y_1 \ast y_2)^{-1},
\end{eqnarray*}
where  $y_1, y_2$ run over all elements of  $X$ satisfying
$y_2 = y_1 \ast (y_2 y_1)$ and $y_2 = y_2 \ast (y_1^2)$.
\end{sect}
{\it Proof.}
Let  $S$  be a surface braid of degree $4$ described by the
$4$-chart $\Gamma$  in
Fig.~\ref{twisttre1holo} (Fig.~\ref{twisttre1cocy}).
Let $\beta_1, \dots, \beta_6$  be the paths from points
in the distinguished regions of the white vertices as in
Fig.~\ref{twisttre1cocy}.
  Let
$W_i$
be the white
vertex near the end of $\beta_i$, for $i=1, \ldots, 6$.
The 4-braids $\rho_S (\beta_1), \dots,
\rho_S (\beta_6)$ are represented by
the intersection braid words
$w_\Gamma (\beta_1), \dots, w_\Gamma (\beta_6)$, which are
$$ \sigma_1, \quad \sigma_2^{-1}\sigma_1, \quad
\sigma_2^{-1}\sigma_1, \quad \sigma_1\sigma_3, \quad
\sigma_1\sigma_3, \quad {\rm and} \quad \sigma_3, $$
respectively.

The quandle automorphisms  $Q(\sigma_1)$, $Q(\sigma_2^{-1}\sigma_1)$,
$Q(\sigma_1\sigma_3)$ and $Q(\sigma_3)$ of
$F_Q \langle x_1, \dots, x_m \rangle$  map the generators as follows.
\begin{eqnarray*}
Q(\sigma_1) &: &
x_1 \mapsto x_2 \ast x_1^{-1}, \quad
x_2 \mapsto x_1, \quad
x_3 \mapsto x_3, \quad
x_4 \mapsto x_4, \\
Q(\sigma_2^{-1}\sigma_1) &: &
x_1 \mapsto x_2 \ast x_1^{-1}, \quad
x_2 \mapsto x_3, \quad
x_3 \mapsto x_1 \ast x_3, \quad
x_4 \mapsto x_4, \\
Q(\sigma_1\sigma_3) &: &
x_1 \mapsto x_2 \ast x_1^{-1}, \quad
x_2 \mapsto x_1, \quad
x_3 \mapsto x_4 \ast x_3^{-1}, \quad
x_4 \mapsto x_3, \\
Q(\sigma_3) &: &
x_1 \mapsto x_1, \quad
x_2 \mapsto x_2, \quad
x_3 \mapsto x_4 \ast x_3^{-1}, \quad
x_4 \mapsto x_3.
\end{eqnarray*}
Recall that the inclusion-induced
quandle homomorphism
$i_\ast :  Q_{y_0}= F_Q \langle x_1, \dots, x_4 \rangle
\to Q(S)$ is the natural projection from
$F_Q \langle x_1, \dots, x_4 \rangle$ to
\begin{eqnarray*}
\langle x_1, \dots, x_4 |
& & x_2 = x_1 \ast (x_2 x_1), \\
& & x_2 = x_2 \ast (x_1^2), \\
& & x_3 = x_2, \\
& & x_4 = x_1
\rangle \\
=
\langle x_1, x_2 |
& & x_2 = x_1 \ast (x_2 x_1), \\
& & x_2 = x_2 \ast (x_1^2)
\rangle .
\end{eqnarray*}

Then the quandle triples of the white vertices  $W_1, \dots,
W_6$  are
\begin{eqnarray*}
(x_2 \ast x_1^{-1}, x_1, x_3) &=&
(x_1 \ast x_2, x_1, x_2), \\
(x_2 \ast x_1^{-1}, x_3, x_1 \ast x_3) &=&
(x_1 \ast x_2, x_2, x_1 \ast x_2), \\
(x_3, x_1 \ast x_3, x_4) &=&
(x_2, x_1 \ast x_2, x_1), \\
(x_1, x_4 \ast x_3^{-1}, x_3) &=&
(x_1, x_1 \ast x_2, x_2), \\
(x_2 \ast x_1^{-1}, x_1, x_4 \ast x_3^{-1}) &=&
(x_1 \ast x_2, x_1, x_1 \ast x_2), \\
(x_1, x_2, x_4 \ast x_3^{-1}) &=&
(x_1, x_2, x_1 \ast x_2),
\end{eqnarray*}
respectively.   The signs of the white vertices are as follows:
$$
\epsilon(W_1) = \epsilon(W_2) = \epsilon(W_3) = +1, \quad
\epsilon(W_4) = \epsilon(W_5) = \epsilon(W_6) = -1.  $$
Therefore we have
\begin{eqnarray*}
\Phi_\theta (\Gamma) = \sum_c
&&
\theta (c(x_1 \ast x_2), c(x_1), c(x_2))
\theta (c(x_1 \ast x_2), c(x_2), c(x_1 \ast x_2)) \\
&&
\theta (c(x_2), c(x_1 \ast x_2), c(x_1))
\theta (c(x_1), c(x_1 \ast x_2), c(x_2))^{-1} \\
&&
\theta (c(x_1 \ast x_2), c(x_1), c(x_1 \ast x_2))^{-1}
\theta (c(x_1), c(x_2), c(x_1 \ast x_2))^{-1},
\end{eqnarray*}
where $c$ runs over all possible quandle homomorphisms
from $Q(S)$  to  $X$.  Hence
\begin{eqnarray*}
\Phi_\theta (\Gamma) = \sum_{y_1,y_2}
&&
\theta (y_1 \ast y_2, y_1, y_2)
\theta (y_1 \ast y_2, y_2, y_1 \ast y_2)
\theta (y_2, y_1 \ast y_2, y_1) \\
&&
\theta (y_1, y_1 \ast y_2, y_2)^{-1}
\theta (y_1 \ast y_2, y_1, y_1 \ast y_2)^{-1}
\theta (y_1, y_2, y_1 \ast y_2)^{-1},
\end{eqnarray*}
where  $y_1, y_2$ run over all elements of  $X$ satisfying
$y_2 = y_1 \ast (y_2 y_1)$ and $y_2 = y_2 \ast (y_1^2)$.

Since the closure of $S$ in ${\bf R}^4$ is ambient isotopic
to the 2-twist spun trefoil,
the result follows. $\Box$

\begin{sect} {\bf Corollary.\/}
\label{firstcor}
Let  $\theta \in Z^3( R_3; {\bf Z}_3 )$
be
the
3-cocycle
$$ \displaystyle{t^{
- \chi_{(0,1,0)} + \chi_{(0,2,0)} -
\chi_{(0,2,1)}
+ \chi_{(1,0,1)} + \chi_{(1,0,2)} + \chi_{(2,0,2)} + \chi_{(2,1,2)}
}},$$
where  $R_3$  is the dihedral quandle of three elements,
${\bf Z}_3$
  is the cyclic group $\langle t | t^3 = 1 \rangle$
of order three, and $\chi_{(i,j,k)}$'s are characteristic functions
as before.   If  $F$ is the 2-twist spun trefoil, then
$$ \Phi_\theta (F) = 3 + 6t \quad \in {\bf Z}[t, t^{-1}]/ (t^3-1). $$
\end{sect}
{\it Proof.\/}
Every pair $\{ y_1, y_2 \}$ of elements of $R_3$ satisfies
the condition of Theorem~\ref{generalform}.
We have the result by
a direct calculation (Table~1 will be helpful). $\Box$

\begin{center}
\begin{tabular}{|l|l|*{6}{c|}c| } \hline
   {$y_1$} & {$y_2$} &
   {$\theta (y_1 \ast y_2, y_1,  $} &
   {$\theta (y_1 \ast y_2, y_2,  $} &
   {$\theta (y_2, y_1 \ast y_2,  $} &
   {$\theta (y_1, y_1 \ast y_2,  $} &
   {$\theta (y_1 \ast y_2, y_1,  $} &
   {$\theta (y_1, y_2, $} &
   {Prod} \\
     & &
   {$ y_2) $} &
   {$ y_1 \ast y_2) $} &
   {$ y_1) $} &
   {$ y_2)^{-1} $} &
   {$ y_1 \ast y_2)^{-1} $} &
   {$ y_1 \ast y_2)^{-1} $} & \\  \hline
$0$ & $0$ &
  $\theta (0,0,0)  $ &
  $\theta (0,0,0)  $ &
  $\theta (0,0,0)  $ &
  $\theta (0,0,0)^{-1}  $ &
  $\theta (0,0,0)^{-1}  $ &
  $\theta (0,0,0)^{-1}  $ & \\
  & &
  $= t^0 $ &
  $= t^0 $ &
  $= t^0 $ &
  $= t^0 $ &
  $= t^0 $ &
  $= t^0 $ &
  $ t^0 $   \\ \hline
$0$ & $1$ &
  $\theta (2,0,1)  $ &
  $\theta (2,1,2)  $ &
  $\theta (1,2,0)  $ &
  $\theta (0,2,1)^{-1}  $ &
  $\theta (2,0,2)^{-1}  $ &
  $\theta (0,1,2)^{-1}  $ & \\
  & &
  $= t^0 $ &
  $= t^1 $ &
  $= t^0 $ &
  $= t^1 $ &
  $= t^{-1} $ &
  $= t^0 $ &
  $ t^1 $   \\ \hline
$0$ & $2$ &
  $\theta (1,0,2)  $ &
  $\theta (1,2,1)  $ &
  $\theta (2,1,0)  $ &
  $\theta (0,1,2)^{-1}  $ &
  $\theta (1,0,1)^{-1}  $ &
  $\theta (0,2,1)^{-1}  $ & \\
  & &
  $= t^1 $ &
  $= t^0 $ &
  $= t^0 $ &
  $= t^0 $ &
  $= t^{-1} $ &
  $= t^1 $ &
  $ t^1 $   \\ \hline
$1$ & $0$ &
  $\theta (2,1,0)  $ &
  $\theta (2,0,2)  $ &
  $\theta (0,2,1)  $ &
  $\theta (1,2,0)^{-1}  $ &
  $\theta (2,1,2)^{-1}  $ &
  $\theta (1,0,2)^{-1}  $ & \\
  & &
  $= t^0 $ &
  $= t^1 $ &
  $= t^{-1} $ &
  $= t^0 $ &
  $= t^{-1} $ &
  $= t^{-1} $ &
  $ t^1 $   \\ \hline
$1$ & $1$ &
  $\theta (1,1,1)  $ &
  $\theta (1,1,1)  $ &
  $\theta (1,1,1)  $ &
  $\theta (1,1,1)^{-1}  $ &
  $\theta (1,1,1)^{-1}  $ &
  $\theta (1,1,1)^{-1}  $ & \\
  & &
  $= t^0 $ &
  $= t^0 $ &
  $= t^0 $ &
  $= t^0 $ &
  $= t^0 $ &
  $= t^0 $ &
  $ t^0 $   \\ \hline
$1$ & $2$ &
  $\theta (0,1,2)  $ &
  $\theta (0,2,0)  $ &
  $\theta (2,0,1)  $ &
  $\theta (1,0,2)^{-1}  $ &
  $\theta (0,1,0)^{-1}  $ &
  $\theta (1,2,0)^{-1}  $ & \\
  & &
  $= t^0 $ &
  $= t^1 $ &
  $= t^0 $ &
  $= t^{-1} $ &
  $= t^1 $ &
  $= t^0 $ &
  $ t^1 $   \\ \hline
$2$ & $0$ &
  $\theta (1,2,0)  $ &
  $\theta (1,0,1)  $ &
  $\theta (0,1,2)  $ &
  $\theta (2,1,0)^{-1}  $ &
  $\theta (1,2,1)^{-1}  $ &
  $\theta (2,0,1)^{-1}  $ & \\
  & &
  $= t^0 $ &
  $= t^1 $ &
  $= t^0 $ &
  $= t^0 $ &
  $= t^0 $ &
  $= t^0 $ &
  $ t^1 $   \\ \hline
$2$ & $1$ &
  $\theta (0,2,1)  $ &
  $\theta (0,1,0)  $ &
  $\theta (1,0,2)  $ &
  $\theta (2,0,1)^{-1}  $ &
  $\theta (0,2,0)^{-1}  $ &
  $\theta (2,1,0)^{-1}  $ & \\
  & &
  $= t^{-1}  $ &
  $= t^{-1} $ &
  $= t^1 $ &
  $= t^0 $ &
  $= t^{-1} $ &
  $= t^0 $ &
  $ t^1 $   \\ \hline
$2$ & $2$ &
  $\theta (2,2,2)  $ &
  $\theta (2,2,2)  $ &
  $\theta (2,2,2)  $ &
  $\theta (2,2,2)^{-1}  $ &
  $\theta (2,2,2)^{-1}  $ &
  $\theta (2,2,2)^{-1}  $ & \\
  & &
  $= t^0 $ &
  $= t^0 $ &
  $= t^0 $ &
  $= t^0 $ &
  $= t^0 $ &
  $= t^0 $ &
  $ t^0 $   \\ \hline

\end{tabular}

\vspace{3mm}

{\bf Table 1 }

\end{center}

\begin{figure}
\begin{center}
\mbox{
\epsfxsize=3in
\epsfbox{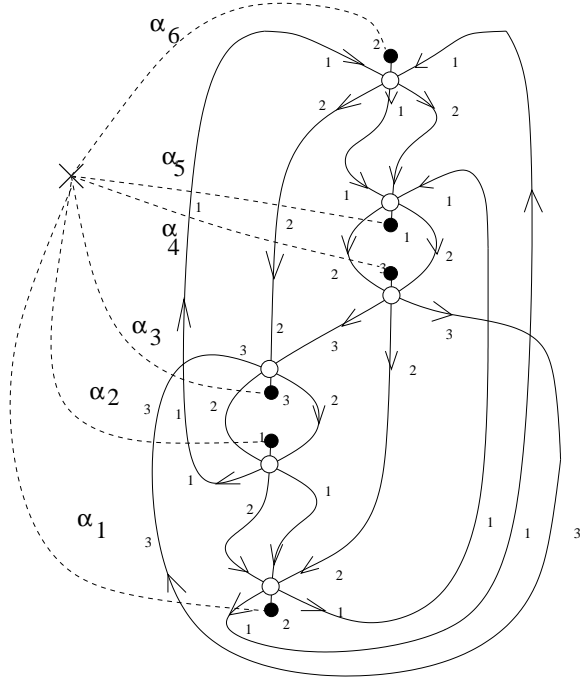}
} \end{center}
\caption{The Hurwitz system for the orientation reversed image
}
\label{twisttre2holo}
\end{figure}

\begin{sect} {\bf Theorem.\/}
\label{othergenform}
Let  $F'$
be the $2$-twist spun trefoil whose
orientation is reversed,
and  $\theta$ a
quandle $3$-cocycle of a finite quandle  $X$ with
coefficient group  $A$.  The
state-sum
invariant
$\Phi_\theta (F')$ is
\begin{eqnarray*}
\sum_{y_1,y_2}
&&
\theta (y_2, y_1 \ast y_2, y_1)^{-1}
\theta (y_2, y_1, y_2)^{-1}
\theta (y_1, y_2, y_1 \ast y_2)^{-1} \\
&&
\theta (y_1 \ast y_2, y_2, y_1)
\theta (y_2, y_1 \ast y_2, y_2)
\theta (y_1 \ast y_2, y_1, y_2),
\end{eqnarray*}
where  $y_1, y_2$ run over all elements of  $X$ satisfying
$y_2 = y_1 \ast (y_2 y_1)$ and $y_2 = y_2 \ast (y_1^2)$.
\end{sect}
{\it Proof.\/}
Let
$S'$
be a surface braid of degree  $4$  described by
a $4$-chart
$\Gamma'$
in Fig.~\ref{twisttre2holo}.  It is known that the closure
of
$S'$
is ambient isotopic to the 2-twist spun trefoil
with the reversed orientation.  (In general, if an
$m$-chart  $\Gamma_2$  is a mirror image of another  $\Gamma_1$,
then the closure of a surface braid described by  $\Gamma_2$
is ambient isotopic to
the closure of a surface braid described by  $\Gamma_1$
whose orientation is reversed.)

For a Hurwitz arc system ${\cal A} = ( \alpha_1, \dots,
\alpha_6 )$ as in Fig.~\ref{twisttre2holo},
the braid system
$$( w_1^{-1} \sigma_{k_1}^{\epsilon_1} w_1, \dots,
w_6^{-1} \sigma_{k_6}^{\epsilon_6} w_6 )$$
of  $S'$
is given as follows:
$$
\begin{array}{ll}

w_1 = \sigma_1\sigma_3^{-1}, \quad &
\sigma_{k_1}^{\epsilon_1} = \sigma_2, \\

w_2 = \sigma_2\sigma_1^{-1}\sigma_3^{-1}, \quad &
\sigma_{k_2}^{\epsilon_2} = \sigma_1^{-1}, \\

w_3 = \sigma_2\sigma_1^{-1}\sigma_3^{-1}, \quad &
\sigma_{k_3}^{\epsilon_3} = \sigma_3, \\

w_4 = \sigma_2^2\sigma_1^{-1}, \quad &
\sigma_{k_4}^{\epsilon_4} = \sigma_3^{-1}, \\

w_5 = \sigma_2^2\sigma_1^{-1}, \quad &
\sigma_{k_5}^{\epsilon_5} = \sigma_1, \\

w_6 = 1, \quad &
\sigma_{k_6}^{\epsilon_6} = \sigma_2^{-1}.

\end{array}
$$

The quandle automorphisms
$Q(\sigma_1\sigma_3^{-1})$,
$Q(\sigma_2\sigma_1^{-1}\sigma_3^{-1})$,
$Q(\sigma_2^2\sigma_1^{-1})$, and
$Q(1)$
of
$F_Q \langle x_1, \dots, x_m \rangle$  map the generators as follows.

\begin{eqnarray*}
Q(\sigma_1\sigma_3^{-1}) &: &
x_1 \mapsto x_2 \ast x_1^{-1}, \quad
x_2 \mapsto x_1, \quad
x_3 \mapsto x_4, \quad
x_4 \mapsto x_3 \ast x_4, \\
Q(\sigma_2
\sigma_1^{-1}\sigma_3^{-1}) &: &
x_1 \mapsto x_2, \quad
x_2 \mapsto x_4 \ast (x_2^{-1} x_1^{-1} x_2), \quad
x_3 \mapsto x_1 \ast x_2, \quad
x_4 \mapsto x_3 \ast x_4, \\
Q(\sigma_2^2\sigma_1^{-1}) &: &
x_1 \mapsto x_2, \quad
x_2 \mapsto x_1 \ast (x_2 x_3^{-1} x_2^{-1} x_1^{-1} x_2), \quad
x_3 \mapsto x_3 \ast (x_2^{-1} x_1^{-1} x_2), \quad
x_4 \mapsto x_4, \\
Q(1)
&: &
x_1 \mapsto x_1, \quad
x_2 \mapsto x_2, \quad
x_3 \mapsto x_3, \quad
x_4 \mapsto x_4.
\end{eqnarray*}

Hence the defining relations
$Q(w_i)( x_{k_i} ) = Q(w_i)( x_{k_i +1} )$ ($i=1,\dots, 6$) of
$Q(S')$
are
\begin{eqnarray*}
x_1  &=&
x_4, \\
x_2  &=&
x_4 \ast (x_2^{-1} x_1^{-1} x_2),  \\
x_1 \ast x_2  &=&
x_3 \ast x_4, \\
x_3 \ast (x_2^{-1} x_1^{-1} x_2)  &=&
x_4, \\
x_2   &=&
x_1 \ast (x_2 x_3^{-1} x_2^{-1} x_1^{-1} x_2), \\
x_2   &=&
x_3.
\end{eqnarray*}

Thus the quandle $Q(S')$ is
\begin{eqnarray*}
\langle x_1, \dots, x_4 |
& & x_2 = x_1 \ast (x_2 x_1), \\
& & x_2 = x_2 \ast (x_1^2), \\
& & x_3 = x_2, \\
& & x_4 = x_1
\rangle \\
=
\langle x_1, x_2 |
& & x_2 = x_1 \ast (x_2 x_1), \\
& & x_2 = x_2 \ast (x_1^2)
\rangle.
\end{eqnarray*}

\begin{figure}
\begin{center}
\mbox{
\epsfxsize=3in
\epsfbox{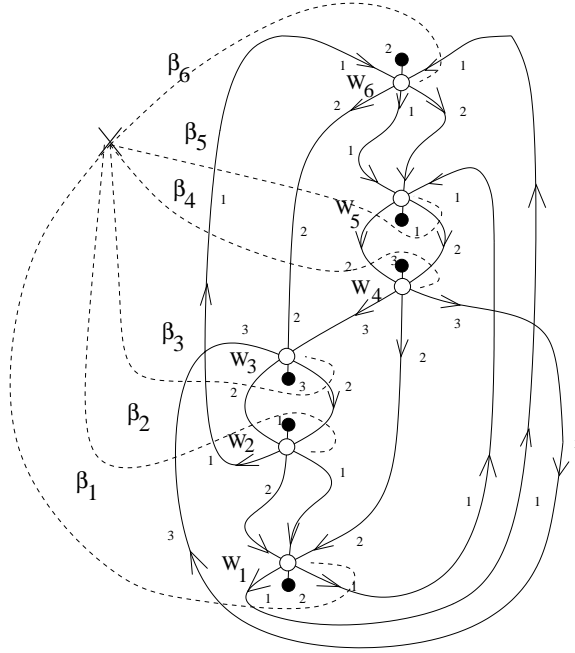}
} \end{center}
\caption{Computing cocycle invariants for the orientation reversed image
}
\label{twisttre2cocy}
\end{figure}

Let $\beta_1, \dots, \beta_6$  be the paths from points
in the distinguished regions of the white vertices as in
Fig.~\ref{twisttre2cocy}, and let
$W_i$ be the white vertex near the end of $\beta_i$ for $i=1, \ldots, 6$.
The 4-braids $\rho_S (\beta_1), \dots,
\rho_S (\beta_6)$, which are represented by
the intersection braid words
$w_\Gamma (\beta_1), \dots, w_\Gamma (\beta_6)$, are
$$
\sigma_1^2\sigma_3^{-1}, \quad
\sigma_2^2\sigma_1^{-1}\sigma_3^{-1}, \quad
\sigma_2^2\sigma_1^{-1}\sigma_3^{-1}, \quad
\sigma_2^3\sigma_1^{-1}, \quad
\sigma_2^3\sigma_1^{-1}, \quad    {\rm and} \quad
\sigma_1,
$$
respectively.
The quandle automorphisms
$Q(\sigma_1^2\sigma_3^{-1})$,
$Q(\sigma_2^2\sigma_1^{-1}\sigma_3^{-1})$,
$Q(\sigma_2^3\sigma_1^{-1})$, and
$Q(\sigma_1)$
of
$F_Q \langle x_1, \dots, x_m \rangle$  map the generators as follows.

\begin{eqnarray*}
Q(\sigma_1^2\sigma_3^{-1}) &: &
x_1 \mapsto x_1 \ast (x_2^{-1} x_1^{-1}), \quad
x_2 \mapsto x_2 \ast x_1^{-1}, \quad
x_3 \mapsto x_4, \quad
x_4 \mapsto x_3 \ast x_4, \\
Q(\sigma_2^2\sigma_1^{-1}\sigma_3^{-1}) &: &
x_1 \mapsto x_2, \quad
x_2 \mapsto x_1 \ast (x_2 x_4^{-1} x_2^{-1} x_1^{-1} x_2), \quad
x_3 \mapsto x_4 \ast (x_2^{-1} x_1^{-1} x_2), \quad
x_4 \mapsto x_3 \ast x_4, \\
Q(\sigma_2^3\sigma_1^{-1}) &: &
x_1 \mapsto x_2, \quad
x_2 \mapsto x_3 \ast (x_2^{-1} x_1^{-1} x_2 x_3^{-1} x_2^{-1} x_1^{-1} x_2),
\\ & &
x_3 \mapsto x_1 \ast (x_2 x_3^{-1} x_2^{-1} x_1^{-1} x_2), \quad
x_4 \mapsto x_4, \\
Q(\sigma_1) &: &
x_1 \mapsto x_2 \ast x_1^{-1}, \quad
x_2 \mapsto x_1, \quad
x_3 \mapsto x_3, \quad
x_4 \mapsto x_4.
\end{eqnarray*}

Then the quandle triples of the white vertices  $W_1, \dots,
W_6$
are
\begin{eqnarray*}
(x_1 \ast (x_2^{-1} x_1^{-1}),
x_2 \ast x_1^{-1},
x_4)
&=&
(x_2, x_1 \ast x_2, x_1), \\
(x_2,
x_1 \ast (x_2 x_4^{-1} x_2^{-1} x_1^{-1} x_2),
x_4 \ast (x_2^{-1} x_1^{-1} x_2))
&=&
(x_2, x_1, x_2), \\
(x_1 \ast (x_2 x_4^{-1} x_2^{-1} x_1^{-1} x_2),
x_4 \ast (x_2^{-1} x_1^{-1} x_2),
x_3 \ast x_4)
&=&
(x_1, x_2, x_1 \ast x_2), \\
(x_3 \ast (x_2^{-1} x_1^{-1} x_2 x_3^{-1} x_2^{-1} x_1^{-1} x_2),
x_1 \ast (x_2 x_3^{-1} x_2^{-1} x_1^{-1} x_2),
x_4)
&=&
(x_1 \ast x_2, x_2, x_1), \\
(x_2,
x_3 \ast (x_2^{-1} x_1^{-1} x_2 x_3^{-1} x_2^{-1} x_1^{-1} x_2),
x_1 \ast (x_2 x_3^{-1} x_2^{-1} x_1^{-1} x_2))
&=&
(x_2, x_1 \ast x_2, x_2), \\
(x_2 \ast x_1^{-1},
x_1, x_3)
&=&
(x_1 \ast x_2, x_1, x_2),
\end{eqnarray*}
respectively.   The signs of the white vertices are as follows:
$$
\epsilon(W_1) = \epsilon(W_2) = \epsilon(W_3) = -1, \quad
\epsilon(W_4) = \epsilon(W_5) = \epsilon(W_6) = +1.  $$
Therefore we have
\begin{eqnarray*}
\Phi_\theta ( \Gamma'
) = \sum_{y_1,y_2}
&&
\theta (y_2, y_1 \ast y_2, y_1)^{-1}
\theta (y_2, y_1, y_2)^{-1}
\theta (y_1, y_2, y_1 \ast y_2)^{-1} \\
&&
\theta (y_1 \ast y_2, y_2, y_1)
\theta (y_2, y_1 \ast y_2, y_2)
\theta (y_1 \ast y_2, y_1, y_2),
\end{eqnarray*}
where  $y_1, y_2$ run over all elements of  $X$ satisfying
$y_2 = y_1 \ast (y_2 y_1)$ and $y_2 = y_2 \ast (y_1^2)$.
This completes the proof. $\Box$

\begin{sect} {\bf Corollary.\/}
\label{secondcor}
Let  $\theta \in Z^3( R_3; {\bf Z}_3 )$  be
the
3-cocycle
$$ \displaystyle{t^{
- \chi_{(0,1,0)} + \chi_{(0,2,0)} -
\chi_{(0,2,1)}
+ \chi_{(1,0,1)} + \chi_{(1,0,2)} + \chi_{(2,0,2)} + \chi_{(2,1,2)}
}},$$
where  $R_3$  is the dihedral quandle of three elements,
${\bf Z}_3$  is the cyclic group $\langle t | t^3 = 1 \rangle$
of order three, and $\chi_{(i,j,k)}$'s are characteristic functions
as before.   If  $F'$
is the 2-twist spun trefoil
whose orientation is reversed, then
$$\Phi_\theta (F')
= 3 + 6t^2 \quad \in {\bf Z}[t, t^{-1}]/(t^3-1).$$
\end{sect}
{\it Proof.\/}
Every pair $\{ y_1, y_2 \}$ of elements of $R_3$ satisfies
the condition of Theorem~\ref{othergenform}.
We have the result by
a direct calculation (Table~2 will be helpful). $\Box$

\vspace{5mm}

\begin{center}
\begin{tabular}{|l|l|*{6}{c|}c| } \hline
   {$y_1$} & {$y_2$} &
   {$\theta (y_2, y_1 \ast y_2,  $} &
   {$\theta (y_2, y_1,  $} &
   {$\theta (y_1, y_2,  $} &
   {$\theta (y_1 \ast y_2, y_2,  $} &
   {$\theta (y_2, y_1 \ast y_2,  $} &
   {$\theta (y_1 \ast y_2, y_1, $} &
   {Prod} \\
     & &
   {$ y_1)^{-1} $} &
   {$ y_2)^{-1} $} &
   {$ y_1 \ast y_2)^{-1} $} &
   {$ y_1 ) $} &
   {$ y_2) $} &
   {$ y_2) $} & \\  \hline
$0$ & $0$ &
  $\theta (0,0,0)^{-1}  $ &
  $\theta (0,0,0)^{-1}  $ &
  $\theta (0,0,0)^{-1}  $ &
  $\theta (0,0,0)   $ &
  $\theta (0,0,0)   $ &
  $\theta (0,0,0)   $ & \\
  & &
  $= t^0 $ &
  $= t^0 $ &
  $= t^0 $ &
  $= t^0 $ &
  $= t^0 $ &
  $= t^0 $ &
  $ t^0 $   \\ \hline
$0$ & $1$ &
  $\theta (1,2,0)^{-1}  $ &
  $\theta (1,0,1)^{-1}  $ &
  $\theta (0,1,2)^{-1}  $ &
  $\theta (2,1,0)   $ &
  $\theta (1,2,1)   $ &
  $\theta (2,0,1)   $ & \\
  & &
  $= t^0 $ &
  $= t^{-1} $ &
  $= t^0 $ &
  $= t^0 $ &
  $= t^0 $ &
  $= t^0 $ &
  $ t^2 $   \\ \hline
$0$ & $2$ &
  $\theta (2,1,0)^{-1}  $ &
  $\theta (2,0,2)^{-1}  $ &
  $\theta (0,2,1)^{-1}  $ &
  $\theta (1,2,0)   $ &
  $\theta (2,1,2)   $ &
  $\theta (1,0,2)   $ & \\
  & &
  $= t^0 $ &
  $= t^{-1} $ &
  $= t^{1} $ &
  $= t^0 $ &
  $= t^1 $ &
  $= t^1 $ &
  $ t^2 $   \\ \hline
$1$ & $0$ &
  $\theta (0,2,1)^{-1}  $ &
  $\theta (0,1,0)^{-1}  $ &
  $\theta (1,0,2)^{-1}  $ &
  $\theta (2,0,1)   $ &
  $\theta (0,2,0)   $ &
  $\theta (2,1,0)   $ & \\
  & &
  $= t^{1} $ &
  $= t^1 $ &
  $= t^{-1} $ &
  $= t^0 $ &
  $= t^1 $ &
  $= t^0 $ &
  $ t^2 $   \\ \hline
$1$ & $1$ &
  $\theta (1,1,1)^{-1}  $ &
  $\theta (1,1,1)^{-1}  $ &
  $\theta (1,1,1)^{-1}  $ &
  $\theta (1,1,1)   $ &
  $\theta (1,1,1)   $ &
  $\theta (1,1,1)   $ & \\
  & &
  $= t^0 $ &
  $= t^0 $ &
  $= t^0 $ &
  $= t^0 $ &
  $= t^0 $ &
  $= t^0 $ &
  $ t^0 $   \\ \hline
$1$ & $2$ &
  $\theta (2,0,1)^{-1}  $ &
  $\theta (2,1,2)^{-1}  $ &
  $\theta (1,2,0)^{-1}  $ &
  $\theta (0,2,1)   $ &
  $\theta (2,0,2)   $ &
  $\theta (0,1,2)   $ & \\
  & &
  $= t^0 $ &
  $= t^{-1} $ &
  $= t^0 $ &
  $= t^{-1} $ &
  $= t^1 $ &
  $= t^0 $ &
  $ t^2 $   \\ \hline
$2$ & $0$ &
  $\theta (0,1,2)^{-1}  $ &
  $\theta (0,2,0)^{-1}  $ &
  $\theta (2,0,1)^{-1}  $ &
  $\theta (1,0,2)   $ &
  $\theta (0,1,0)   $ &
  $\theta (1,2,0)   $ & \\
  & &
  $= t^0 $ &
  $= t^{-1} $ &
  $= t^0 $ &
  $= t^1 $ &
  $= t^{-1} $ &
  $= t^0 $ &
  $ t^2 $   \\ \hline
$2$ & $1$ &
  $\theta (1,0,2)^{-1}  $ &
  $\theta (1,2,1)^{-1}  $ &
  $\theta (2,1,0)^{-1}  $ &
  $\theta (0,1,2)   $ &
  $\theta (1,0,1)   $ &
  $\theta (0,2,1)   $ & \\
  & &
  $= t^{-1} $ &
  $= t^0 $ &
  $= t^0 $ &
  $= t^0 $ &
  $= t^1 $ &
  $= t^{-1}  $ &
  $ t^2 $   \\ \hline
$2$ & $2$ &
  $\theta (2,2,2)^{-1}  $ &
  $\theta (2,2,2)^{-1}  $ &
  $\theta (2,2,2)^{-1}  $ &
  $\theta (2,2,2)   $ &
  $\theta (2,2,2)   $ &
  $\theta (2,2,2)   $ & \\
  & &
  $= t^0 $ &
  $= t^0 $ &
  $= t^0 $ &
  $= t^0 $ &
  $= t^0 $ &
  $= t^0 $ &
  $ t^0 $   \\ \hline

\end{tabular}

\vspace{3mm}

{\bf Table 2}

\end{center}

The
Corollaries~\ref{firstcor} and \ref{secondcor}
imply

\begin{sect} {\bf Theorem.\/}
The $2$-twist spun trefoil is non-invertible.
\quad $\Box$  
\end{sect}

Examples 10 and 11 of \cite{FoxTrip}
are inverses of each other.
One has its Alexander ideal generated by $2T-1$; the other has $T-2$ as the
generator of the Alexander ideal. So
the non-invertibility of this (Example 10/11) ribbon knot is detected
by the Alexander ideal. The knot quandle, which contains the fundamental
group and a choice of positive meridional element,
can be used to compute the Alexander ideal. In the
case of the $2$-twist spun trefoil
(which happens to be
Example 12 of \cite{FoxTrip}),
we have computed
that the knotted sphere and its
orientation reversed copy
have the same
knot quandles --- thus they have the same Alexander ideal which is
(non-principally) generated by
$2T-1$ and $T-2.$ The invariant $\Phi_\theta$ is the first known
state-sum
invariant that detects
non-invertibility of this important example.

\end{document}